\documentclass[11pt]{article}
\title{On the Localisation Theorem for rational Cherednik algebra modules}
\author{Rollo Jenkins}
\date{20/01/2014}

\usepackage[a4paper,margin=1in]{geometry}
\usepackage{amsmath}
\usepackage{amsthm}
\usepackage{amssymb}
\usepackage{mathrsfs}
\usepackage{manfnt}
\usepackage{mathtools}
\usepackage{booktabs}    
\usepackage{graphicx}
\usepackage{pxfonts}
\usepackage{stmaryrd}
\usepackage{lscape}
\usepackage[table,dvipsnames]{xcolor}
\definecolor{beige1}{RGB}{253,245,230}
\usepackage[boxsize=5pt]{ytableau}
\usepackage{multirow}
\usepackage[toc,page]{appendix}
\usepackage{enumerate}

\usepackage{url}
\makeatletter
\def\url@leostyle{%
  \@ifundefined{selectfont}{\def\UrlFont{\sf}}{\def\UrlFont{\small\ttfamily}}}
\makeatother 
\usepackage[all]{xy}
\SelectTips{cm}{}
\newdir^{c}{{}*!/-5pt/@^{(}}

\newtheorem{defi}{Definition}[section]
\newtheorem{theo}{Theorem}[section]
\newtheorem{prop}{Proposition}[section]
\newtheorem{exam}{Example}[section]
\newtheorem{coro}{Corollary}[section]
\newtheorem{lemm}{Lemma}[section]
\newtheorem{rema}{Remark}[section]

\newtheorem*{conv}{Convention}
\newtheorem*{introtheo}{Theorem}

\newcommand{\Z}{\mathbb Z}

\newcommand{\R}{\mathbb R}
\newcommand{\CC}{\mathbb C}

\newcommand{\N}{\mathbb N}

\newcommand{\h}{\mathfrak h}

\newcommand{\bydef}{\mathrel{\mathop:}=}

\newcommand{\cc}{\mathbf c}

\newcommand{\fr}[1]{\tfrac{1}{#1}}

\newcommand{\wt}{\mathbf{wt}}

\newcommand{\us}{{X^{\mathrm{us}}}}
\newcommand{\git}[1]{/\!\!/_{\mbox{\!\tiny{${#1}$}}}}

\colorlet{cite}{LimeGreen!50!Green}
\usepackage[%
	bookmarks=true,			
	unicode=true,			
	pdftitle={},		%
	pdfauthor={},	%
	pdfkeywords={} {},	
	colorlinks=true,		
	linkcolor=Red,			
	citecolor=cite,		
	filecolor=magenta,		
	urlcolor=RoyalBlue,			
  ]{hyperref}

\begin{document}

\flushbottom
 \maketitle

\abstract{
Let $W$ be a complex reflection group of the form $G(l,1,n)$. Following \cite{BK, BPW, GordonRemark, GS1, GS2, KR, MN}, the theory of deform quantising conical symplectic resolutions allows one to study the category of modules for the spherical Cherednik algebra, $U_\cc(W)$, via a functor, $\mathbb T_{\cc,\theta}$, which takes invariant global sections of certain twisted sheaves on some Nakajima quiver variety $Y_\theta$. 

A parameter for the Cherednik algebra, $\cc$, is considered \textbf{good} if there exists a choice of GIT parameter $\theta$, such that $\mathbb T_{\cc,\theta}$ is exact and \textbf{bad} otherwise. By calculating the Kirwan--Ness strata for $\theta=\pm(1,\ldots,1)$ and using criteria of \cite{MN2}, it is shown that the set of all bad parameters is bounded. The criteria are then used to show that, for the cases $W=\mathfrak S_n, \mu_3, B_2$, all parameters are good.
}

\tableofcontents
\newpage

\section{Introduction}

Fix $l,n\in \N$, not both one, and define the complex reflection group $W=G(l,1,n)=\mu_l \wr \mathfrak S_n$, the $n$-fold wreath product of the cyclic group of order $l$. Also fix, for the remainder of the paper, $F=\CC^*$. Throughout, the term `rational Cherednik algebra' is understood to mean `rational Cherednik algebra with $t=1$. 

The motivation for this work is to provide a means to study the representation theory of the rational Cherednik algebra, $H_\cc(W)$, by constructing a category of geometric objects. In fact, the theory produces modules for the \emph{spherical subalgebra}, $U_\cc(W)$ of $H_\cc(W)$; see Section \ref{sect:rca}. The parameter $\cc$ is called \textbf{spherical} (or sometimes in the literature \emph{regular}) if the spherical subalgebra has finite global dimension. In such a case, $U_\cc(W)$ and $H_\cc(W)$ are Morita equivalent; see \cite[Theorem 5.5]{Etingof2}.

With the data $(l,n)$, define a quiver, $Q^l_\infty$, and a dimension vector, $\epsilon$; see Section \ref{G::sect:nakajima}. This gives a smooth symplectic variety, $X\bydef T^* \mathrm{Rep}(Q^l_\infty,\epsilon)$, the cotangent bundle on the affine space of all complex representations of $Q^l_\infty$ with dimension vector $\epsilon$. This comes with the action (basis change) of the reductive group $G=\prod_{i=0}^{l-1}\mathrm{GL}_n(\CC))$. Let $F=\CC^*$ act on the coordinate functions of $X$ with weight $1$ and choose a character $\theta$ of $G$ so that the corresponding GIT quotient, $Y_\theta \bydef X\git{\theta} G$ is smooth. Since the $G$ and $F$ actions on $X$ commute, an $F$-action is induced on the quotient $Y_\theta$. 

Let $\chi\in (\mathfrak g^*)^G$ be a character of $\mathfrak g$, this determines a parameter $\cc$ for the rational Cherednik algebra; see Theorem \ref{G::theo:chi and params}. The sheaf of regular functions on $Y_\theta$ is deformed to produce a so-called \emph{W-algebra}, $\mathcal W_{Y_\theta}$, on $Y_\theta$; see Equation \ref{W algebra on Y}. This construction depends on the choice of $\chi$. Finally, consider the category $\left(\widetilde{\mathcal W}_{Y_\theta},F\right)^\mathsf{good}$ of good, $F$-equivariant sheaves of $\CC(\hbar^{1/2})\otimes_{\CC(\hbar)}\mathcal W_{Y_\theta}$-modules; these are defined in \cite[Section 2]{KR}. This category is the candidate geometric analogue of finitely generated $U_\cc(W)$-modules. Indeed, given any $\mathcal M\in \left(\widetilde{\mathcal W}_{Y_\theta},F\right)^\mathsf{good}$ there is an action of $U_\cc(W)$ on the $F$-invariant global sections of $\mathcal M$. This gives a functor,
\[
\xymatrix@C=1in{
\left(\widetilde{\mathcal W}_{Y_\theta},F\right)^\mathsf{good} \ar[r]^{\mathbb T\bydef \Gamma(Y_\theta, -)^F} & {}_{U_\cc(W)}\mathbf{mod}}
\] which depends on both the parameters $\theta$ and $\chi$.

In \cite{MN2}, McGerty and Nevins consider the Kirwan--Ness stratification of the unstable locus of $X$. This is a finite stratification and depends on the choice of GIT parameter $\theta$. To each of these strata they associate a cone of parameter values. If $\chi$ lies outside of all of these cones then the functor $\mathbb T$ is exact for this choice of $\theta$; see Theorem \ref{G::theo:MNT}. 
When the parameter $\chi$ is such that the corresponding rational Cherednik parameter $\cc$ is spherical, this implies that the functor $\mathbb T$ is an abelian equivalence of categories; see Theorem \ref{G::theo:Commutes}.

Say that a fixed choice of parameter $\chi$ is \textbf{bad} if there exists no choice of $\theta$ such that the McGerty--Nevins criteria are satisfied. Theorem \ref{G::theo:essential one-ps} classifies the Kirwan--Ness stratification for 
$\theta=\pm(1,\ldots,1)$, giving a bound on the set of bad parameters. The main result of the paper is a calculation of the Kirwan--Ness strata for $\mu_3$ and $B_2$; applying the criterion then gives the following result.
\begin{introtheo}(Theorems: \ref{G::theo:bad.params.for. S_n},\ \ref{G::theo:bad.params.for.mu3},\ \ref{G::theo:B_2 params}) When $W=\mathfrak S_n$, $W=\mu_3$ or $W=B_2$, for any character $\chi\in (\mathfrak g^* )^G$ there is some $\theta$, not lying on a GIT wall, such that $\mathbb T$ is exact.
\end{introtheo}

Given these calculations, it seems reasonable to conjecture that, for general $l$ and $n$, there exists enough flexibility in choosing $\theta$ so that $\chi$ is never bad in the sense above.

\begin{rema}
This conjecture has recently been proved by Ivan Losev; see \cite[Corollary 2.2]{Losev3}.
\end{rema}

 \section{Background}\label{Sect:background}

\subsection{Rational Cherednik Algebras}\label{sect:rca}

Fix $l,n\in \N$, not both equal to one. Let $W=G(l,1,n)=\mathfrak \mu_l\wr \mathfrak S_n$, with its natural representation $\h$. Let $\mathcal S$ denote the set of non-trivial reflections in $W$ and let $\mathcal E\bydef\{ \ker(1-s)\subset \h\, |\, s\in\mathcal S\}$. Let $\zeta = \exp(\frac{2\pi \sqrt{-1}}{l})$, a primitive $l^{\mathrm{th}}$ root of unity. The group $W$ acts on $\mathcal E$ and on $\CC W$ by conjugation. With respect to these actions, let $\gamma_{(-)}\colon \mathcal E\longrightarrow \CC W;\quad H\mapsto \gamma_H$ be a $W$-equivariant map such that, for each $H\in\mathcal E$, $\gamma_H\in \CC W_H$ and the trace of $\gamma_H$, acting on the $\CC W$-module $\CC W_H$, is zero.

\begin{defi} The \textbf{Rational Cherednik algebra}, $H(W,\h,\gamma)$, associated to the data $(W,\h,\gamma)$ is defined to be the quotient of the smash product of the group, $W$, with the tensor algebra 
$T_\CC(\h\oplus \h^* )$ by the relations
\begin{align*}
 [x,x']&=[y,y']=0 \\ [y,x]&=x(y) +  \sum_{H\in \mathcal E}\frac{\alpha_H(y)x(v_H)}{\alpha_H(v_H)}\gamma_H & &\textrm{for all }y,y'\in\h\textrm{ and }x,x'\in\h^* .
\end{align*}
\end{defi}

There are two ways to parametrise the map $\gamma_{(-)}$ which defines the Cherednik algebra. The first is by parameters, $\cc\bydef (c_0,\ldots,c_{l-1})$, indexed by conjugacy classes of reflections in $\mathcal S$. Following \cite[Section 1.4.1]{Vale}, with the notation defined there, these are defined as 
\begin{align*}
c_0 &\bydef c_{\sigma_{ij}^t} & &\textrm{for all }i,j,t,\\
c_t &\bydef c_{s_k^t} &&\textrm{for all }1\le t\le l-1\textrm{ and all }k.
\end{align*}
These are related to the parameters, $(k,c_{\zeta^1},\ldots,c_{\zeta^{l-1}})$, used in \cite[Section 3.3]{GordonRemark} by, $c_0=-k$ and $c_t=-\fr2 c_{\zeta^t}$ for $t=1,\ldots,l-1$. The corresponding Cherednik algebra is written $H_\cc(W)$.

The second way to parametrise $\gamma_{(-)}$, that will be only be used in Section \ref{Sect:mu3}, is by parameters, $\mathbf k\bydef (k_{00},k_{1},\ldots,k_{l-1})$, indexed by reflection hyperplanes in $\mathcal E$. These are also defined in \cite[Section 1.4.1]{Vale} and are related to the parameters $\cc$ by $k_{00}=-c_0$ and for all $t=1,\ldots,l-1$, $\sum_{j=0}^{l-1}(k_{j+1}-k_j)\zeta^{jt}=-2c_t$, where $k_0=k_{l}=0$. The corresponding Cherednik algebra is written $H_\mathbf k(W)$.

Let $e\bydef\frac{1}{|W|}\sum_{w\in W}w$ be the trivial idempotent of $W$. The \textbf{spherical Cherednik algebra} is the subalgebra $U_\cc(W)\bydef eH_\cc(W) e$ of $H_\cc$. 
The parameter $\cc$ is called \textbf{spherical} if $H_\cc e H_\cc = H_\cc$ and \textbf{aspherical} otherwise. By a theorem of R. Bezrukavnikov, presented in \cite[Theorem 5.5]{Etingof2}, the spherical subalgebra, $U_\cc(W)$, is Morita equivalent to $H_\cc(W)$ precisely when $\cc$ is spherical. In \cite[Theorem 3.3]{DG}, Dunkl and Griffeth give a complete characterisation of aspherical values of $\cc$.

 \subsection{A Nakajima Quiver Variety}\label{G::sect:nakajima}
 Let $Q^l$ and $Q^l_\infty$ denote the following quivers.
\[ 
Q^l\bydef \parbox{1ex}{\xymatrix@R=1.5ex@C=2.5ex{
                                                         & 1 \ar@/^/[dr]^{\mathbf X^{(1)}} & \\
  0 \ar@/^/[ur]^{\mathbf X^{(0)}} &                     & \parbox{1ex}{\vspace{-1ex}\hspace{-.5ex}\vdots} \ar@/^/[dl]^{\mathbf X^{(l-2)}}\\
                                                         & l-1 \ar@/^/[ul]^{\mathbf X^{(l-1)}} &
 }} \quad 
 Q^l_\infty\bydef \parbox{1ex}{\xymatrix@R=1.5ex@C=2.5ex{
                  &                                       & 1 \ar@/^/[dr]^{\mathbf X^{(1)}} & \\
 \infty \ar[r]^v & 0 \ar@/^/[ur]^{\mathbf X^{(0)}} &                     & \parbox{1ex}{\vspace{-1ex}\hspace{-.5ex}\vdots} \ar@/^/[dl]^{\mathbf X^{(l-2)}}\\
                  &                                       & l-1 \ar@/^/[ul]^{\mathbf X^{(l-1)}} &
 }}
 \] Let $\gamma\bydef(n,\ldots,n)$ be a dimension vector for $Q^l$ and $\epsilon=(1,n,\ldots,n)$ for $Q^l_\infty$ (where the first entry corresponds to $\infty$).  
 Let $\overline{Q^l_\infty}$ denote the doubled quiver, with $\mathbf Y^{(i)}$ the reverse of $\mathbf X^{(i)}$ and $w$ the reverse of $v$. Define $V\bydef\mathrm{Rep}(Q_\infty^l,\epsilon)$ and $X\bydef T^* V\cong\mathrm{Rep}\left(\overline{Q^l_\infty}, \epsilon_n\right)$. Note that $X$ is symplectic with respect to the form $\omega_X=\sum_{i,j,k}  d\mathbf X_{ij}^{(k)}\wedge d\mathbf Y_{ij}^{(k)} + \sum_i dv_i\wedge dw_i$. Let $F=\CC^*$ act on $X$ so that each of the coordinate functions has degree one. 

Define $G\bydef \mathrm{GL}_n(\CC)^l$, acting by base change on the representations in $V$. This action is hamiltonian, commutes with the action of $F$ and gives a moment map 
 \[
 \mu\colon X \longrightarrow \mathfrak g \cong \mathfrak g^* ;\qquad \left( \mathbf X,\mathbf Y;v,w\right)\mapsto [\mathbf X,\mathbf Y]+\mathbf v\mathbf w.
 \] See \cite[Section 3.3]{Gordon}. The condition that $\left( \mathbf X,\mathbf Y;v,w\right)\in\mu^{-1}(0)$ is known as the \textbf{ADHM equation}.

Let $\mathcal L$ be the trivial line bundle on $X$ with a choice of $G$-linearisation corresponding to a character $\theta\colon G\longrightarrow \CC^* $. 
That is, an element, $g\in G$, acts on a point, $(x,l)$, in the total space, $X\times \mathbb A^1$, by the rule
\[
g\cdot (x,l)\bydef (g\cdot x,\theta(g)l)\in X\times \mathbb A^1.
\] Let $X\git \theta G\bydef \mathrm{Proj}\left( \oplus_{i\ge 0} \Gamma(X,\mathcal L^{\otimes i})^G\right)$ denote the GIT quotient with respect to this choice of linearisation. The \textbf{Nakajima quiver variety} is the 
 GIT quotient 
 \[
 Y_\theta\bydef \mu^{-1}(0)\git{\theta}G
 \] and the corresponding projective map $p\colon Y_\theta \longrightarrow Y_0$ is a symplectic resolution of singularities when $\theta$ doesn't lie on a GIT wall (see \cite[Section 3.9]{Gordon}). With the action of $F$ introduced above, this is an example of a \emph{conical symplectic resolution}. The points of $Y_\theta$ parametrise equivalence classes of polystable quiver representations of $\overline Q^l_\infty$ that satisfy the ADHM equation.
 The GIT walls in this case have been calculated by Gordon; see \cite[Lemma 4.3 and Remark 4.4]{Gordon}.

Let $\theta=(\theta^0,\ldots,\theta^{l-1})\in \Z^l$ be a GIT parameter. Extend this to a vector 
$\hat \theta=(\theta^\infty,\theta^0,\ldots,\theta^{l-1})$
 so that the dot product, $\hat \theta \cdot \epsilon=0$; that is, choose $\theta^\infty= -\sum_{i=0}^{l-1}n\theta^i$. Given a quiver representation $V$ of $\overline Q^{l}_\infty$, let $\mathbf{dim}V$ be its dimension vector. A proper subrepresentation, $V$, of $x\in X$ is said to \textbf{destabilise} $x$ if $\hat \theta \cdot \mathbf{dim}V<0$. A theorem, \cite[Proposition 3.1]{King}, of King shows that a point, $x\in X$, is semistable with respect to $\theta$ if and only if there does not exist a proper destabilising subrepresentation.

\subsection{Module categories of W-algebras}

Choose a character of $\mathfrak g\bydef \mathrm{Lie}(G)$, $\chi\in (\mathfrak g^* )^G$. The reader is referred to \cite{KR} and \cite{BK} for an introduction to W-algebras on $X$. Let $\mathcal W_X$ denote a W-algebra on $X$ with a formal parameter $\hbar$. Let $\mathbf k\bydef \CC(\hbar)$ and extend the scalars of $\mathcal W_X$ to the field $\mathbf k(\hbar^{1/2})$ by defining
\[
 \widetilde{\mathcal W_X}\bydef \mathbf k(\hbar^{1/2})\otimes_{\mathbf k}\mathcal W_X.
 \]

Let $(\widetilde{\mathcal W_X},F)$ denote the category of $F$-equivariant $\widetilde{\mathcal W_X}$-modules; $(\widetilde{\mathcal W_X},F)^\mathrm{good}$, the full subcategory of good modules. Let 
$(\widetilde{\mathcal W_X},G,F)^\mathrm{good}$ denote the category of good, $F$-equivariant, quasi-$G$-equivariant, $\widetilde{\mathcal W_X}$-modules and $(\widetilde{\mathcal W_X},G,F)^\mathrm{good}_\chi$, its full subcategory of $\chi$-twisted modules. See  \cite[Section 2]{KR} for definitions of these categories.

 \section{Localisation for $U_\cc(W)$}\label{G::chap:localisation.for.Gl1n}
 
Fix a parameter $\cc=(c_0,\ldots,c_{l-1})$\footnote{When $n=1$, the parameter takes the form, $\cc=(c_1,\ldots,c_{l-1})$.} for $H_\cc(W)$.

\subsection{Deformation Quantisations}

The anti-isomorphism $(-)^\mathrm{op}\colon D_V\longrightarrow D_{V^* }$, that maps $\mathbf X_{ij}^{(m)}\mapsto \tfrac{\partial}{\partial \mathbf Y_{ij}^{(m)}}$ and $\tfrac{\partial}{\partial \mathbf X_{ij}^{(m)}}\mapsto \mathbf Y_{ij}^{(m)}$, produces an $F$-equivariant isomorphism, $\phi\colon \mathcal W_X \stackrel{\simeq}{\longrightarrow} \mathcal W_X^\mathrm{op}$, that acts on $\Gamma(X,\mathcal W_X(1))$ by exchanging $\mathbf X_{ij}^{(m)}$ and $\mathbf Y_{ij}^{(m)}$ for all $i,j,m$.

If $a\in \mathcal W_X$ and $g\in G$ then $g\cdot \phi(a)=\phi(g^{-1}\cdot a)$, so $\phi$ restricts to an isomorphism $\mathcal W_X^G \cong (\mathcal W_X^\mathrm{op})^G$. Fix some choice of quantised moment map, $\tau$, and define the algebra
 \[
 D^G_{\tau,\chi}(W)\bydef  \Gamma\left(\frac{\mathcal W_X}{\mathcal W_X\langle \tau(A) - \chi(A)\, |\, A\in \mathfrak g\rangle}\right)^{F,G}\cong  \left(\frac{D_V}{D_V\langle \tau(A) - \chi(A)\, |\, A\in \mathfrak g\rangle}\right)^{G}.
 \]
Deformation quantisations of $Y_\theta$ are of the form 
\[
\mathcal W^\beta \bydef \left( \frac{\mathcal W_X}{\mathcal W_X\langle \tau - \beta\rangle}\right)^G,
\] where $\beta\in (\mathfrak g^* )^G$ is some character. By \cite[Corollary 2.3.3]{Losev}, these are in bijection with $H^2(Y_\theta,\CC)$ via the period map. Under this bijection, if some quantisation, $\mathcal W$ say, is mapped to $\alpha\in H^2(Y_\theta,\CC)$ then $\mathcal W^\mathrm{op}$ is mapped to $-\alpha$. Let the
cohomology class corresponding to $\mathcal W^\beta$ be denoted $\mathrm{Per}(\mathcal W^\beta)$. 
\begin{conv}
Adopt the convention that $\beta$ is chosen so that $\phi(\tau-\beta)=\tau-\beta$ and call the corresponding quantised moment map
\[
\hat\tau \bydef \tau - \beta.
\] 
\end{conv}
This choice of $\beta$ implies that $\phi(\mathcal W_X\hat\tau)=\hat\tau\mathcal W_X$, which in turn implies that
\[
(\mathcal W^\beta)^\mathrm{op}=
\phi\left(\left(\frac{\mathcal W_X}{\mathcal W_X\langle \hat\tau\rangle }\right)^G\right) \cong \left(\frac{\mathcal W_X}{\mathcal W_X\langle \hat\tau\rangle }\right)^G = \mathcal W^\beta.
\] Therefore, $\mathrm{Per}(\mathcal W^\beta)=-\mathrm{Per}(\mathcal W^\beta)$, so that both are zero in $H^2(Y_\theta,\CC)$.
 
 \begin{exam}\label{G::exam:quantised.moment.map.for.vs} Consider $G=\CC^\times$ acting on an $n$-dimensional vector space, $V$, with weights $a_1,\ldots,a_n$ on an eigenbasis $x_1,\ldots,x_n$. Extend this to a hamiltonian action on $X=T^* V$ with $X_1,\ldots,X_n, Y_1,\ldots,Y_n$ coordinate functions as above. Then an arbitrary quantised moment map is of the form $\tau_a(1)=-\sum_ia_iX_iY_i\hbar^{-1}+a\mathrm{tr}(1)=-\sum_ia_iX_iY_i\hbar^{-1}+an$ for some constant $a\in\CC$. Now, $-\fr{2}\sum_ia_i(X_iY_i\hbar^{-1}+Y_iX_i\hbar^{-1})$
 is invariant under $\phi$, so set it to equal $\tau(1)$. Thus, $\hat\tau=-\sum_ia_iX_iY_i\hbar^{-1}-\frac{1}{2}\sum_ia_i= \tau_{-\frac{1}{2n}\sum_ia_i}$.
 \end{exam}

\subsection{The relationship between the parameter $\cc$ for $U_\cc(W)$ and the character $\chi\in (\mathfrak g^* )^G$}

It is now necessary to calculate the relationship between the character, $\chi$, used to twist modules in $D^G_{\tau,\chi}(W)$ and the parameter $\cc$ for the corresponding spherical Cherednik algebra $U_\cc(W)$. 
Let $\zeta$ be a primitive $l^\textrm{th}$ root of unity. Let $I_n$ denote the identity matrix in $\mathrm{GL}_n(\CC)$.
For $i=0,\ldots,l-1$ let $I^{(i)}\bydef(0,\ldots,0,I_n,0,\ldots,0)\in \mathfrak g$. Define the characters of $\mathfrak g$ by, 
 \[
 \mathrm{tr}^{(i)}(A_0,\ldots,A_{l-1})\bydef \mathrm{tr}(A_i)
 \] for $(A_0,\ldots,A_{l-1})\in \mathfrak g$. Let $\{X_{s,t}^{(i)}\, |\, s,t=1,\ldots,n,\, i=0,\ldots,l-1\}$ be the differential operators which multiply by coordinate functions on $V$ and let $\partial_{s,t}^{(i)}$ be the corresponding partial derivatives.

Given a parameter $\cc=(c_0,\ldots,c_{l-1})$, define a character, $\chi\in (\mathfrak g^* )^{G}$, by
 \begin{equation}\label{G::equa:conversion}
 \chi \bydef \begin{cases} \left( c_0+ \frac{1}{2}\right) \mathrm{tr}^{(0)} &\textrm{$l=1$,}\\[1ex]
 \frac{1}{l}\sum_{t=1}^{l-1}\left(1-2\sum_{k=1}^{l-1}\zeta^{kt}c_k\right)\mathrm{tr}^{(t)} &\textrm{$n=1$,} \\[1ex]
  \left(c_0+\frac{1}{2} + \frac{1}{l}\left(-1+l-2\sum_{i=1}^{l-1}c_i\right)\right)\mathrm{tr}^{(0)} + \frac{1}{l}\sum_{t=1}^{l-1}\left(-1-2\sum_{k=1}^{l-1}\zeta^{kt}c_k\right)\mathrm{tr}^{(t)}
 &\parbox{1.2cm}{\textrm{$l>1$,}\\ \textrm{$n>1$.}}
 \end{cases}
 \end{equation}

 The following theorem is a combination of the results of \cite{GordonRemark}, \cite{BK} and \cite{GGS}. 
 \begin{theo}(Gordon, Bellamy--Kuwabara, Ginzburg--Gordon--Stafford)\label{G::theo:chi and params} There is an isomorphism
 \[
 D^G_{\hat\tau,\chi}(W)\cong U_\cc(W).
 \]
 \begin{proof}

  Suppose $l=1$, so that $G=\mathrm{GL}_n(\CC)$ and $V=\mathfrak g\times U$, where $U=\CC^n$. Let $x_1,\ldots,x_n$ be coordinate functions on $U$ and, for each $i=1,\ldots,n$, let $\partial_i$ be the partial derivative with respect to $x_i$. In \cite{GGS}, they choose a quantised moment map, $\tau$, corresponding to the action of $\mathfrak g$ on $\CC[V]$ via derivations. The identity matrix $I_n$ acts on $\CC[V]$ by the derivation $-\sum_ix_i\partial_i$; so that $\tau(I_n)=-\sum_ix_i\partial_i$.
 
 For $g\in G$ and $(M;v)\in V=\mathrm{Mat}_n(\CC)\times \CC^n$, $g\cdot (M;v)=(gMg^{-1};gv)$. Therefore, if $g$ lies in the centre of $G$ it acts trivially on $\mathrm{Mat}_n(\CC)\times \{0\}\subset V$. It follows that the differential of the action of any scalar matrix acts by zero on $\mathrm{Mat}_n(\CC)\times\{0\}$ and so the corresponding vector field given by $\hat\tau$ along this subvariety must be zero. For this reason, Example \ref{G::exam:quantised.moment.map.for.vs} in the special case $a_1=\cdots=a_n=1$ must agree with $\hat\tau(I_n)$. Thus 
 \[
 \hat\tau(I_n)=\hat\tau_{-\frac{1}{2n}\sum_i 1}(1)=-\sum_ix_i\partial_i-\frac{1}{2}\mathrm{tr}(I_n)=\tau(I_n)-\frac{1}{2}\mathrm{tr}(I_n).
 \]
Because $\hat\tau$ and $\tau-\fr2\mathrm{tr}$ agree at $I_n\in\mathfrak g$ and differ by a character they must be equal. Now \cite[Theorem 2.8]{GGS} gives
 \[ 
 D^G_{\hat\tau,\chi} = \left(\frac{D_V}{D_V\langle\hat\tau (A) - \chi(A)\, |\, A\in\mathfrak g\rangle}\right)^G= \left(\frac{D_V}{D_V\langle\tau(A) - (\chi + \fr2)(A)\rangle}\right)^G\cong D^G_{\tau,\chi+\fr2}\cong U_{\chi-\frac{1}{2}}.
 \] Therefore, setting $\chi = c_0 + \fr2$ gives the required result.

  Suppose $n=1$. In this case the first component of the character doesn't contribute to twisting $\mathcal W$-modules in the sense that if $\chi$ and $\chi'$ differ by $\mathrm{tr}^{(0)}$ then $D^G_{\hat\tau,\chi}(W)=D^G_{\hat\tau,\chi'}(W)$. 
  Using the convention $\hat\tau^\mathrm{op}=\hat\tau$, the quantised moment map is
 \begin{align*}
 \hat\tau(I^{(j)})&= \fr 2(X^{(j)}\partial^{(j)}+\partial^{(j)}X^{(j)}) - \fr2(X^{(j-1)}\partial^{(j-1)}+\partial^{(j-1)}X^{(j-1)}) - \fr2(v\partial_v + \partial_vv)\delta_{j,0}\\
   			  &=  	  X^{(j)}\partial^{(j)}  - X^{(j-1)}\partial^{(j-1)} - (v\partial_v +\fr2)\delta_{j,0},
\end{align*} where $j=0,\ldots,l-1$.
 Let $\chi= \sum \chi_i \mathrm{tr}^{(i)}$, an arbitrary character. Then summing $\hat\tau(I^{(j)})-\chi_j$ over all $j$ gives 
 \[
 v\partial_v+\sum_{i=0}^{l-1} \chi_i +\fr2 \in \langle \hat\tau(A) - \chi(A)\, |\, A\in \mathfrak g\rangle.
 \]
 For $i=0,\ldots, l-1$, let $B_i\bydef X^{(i)}\partial^{(i)}$ and let $C\bydef X^{(0)}\cdots X^{(l-1)}$ and $D\bydef \partial^{(0)}\cdots\partial^{(l-1)}$. Then $B_0,\ldots,B_{l-1},C,D$ and $v\partial_v$ generate
 $D_V^G$ and they satisfy the relation $B_0\cdots B_{l-1} = CD$. Therefore, $B_0,\ldots,B_{l-1},C,D$ generate $U_\chi$ and the relations are
\begin{align*} B_0\cdots B_{l-1} &= CD,&  B_1-B_0 &=\chi_1, & &\ldots,& B_{l-1}-B_{l-2}&=\chi_{l-1}.
 \end{align*} 

Now consider another algebra, $A_\chi$, constructed in a different way. Define $V'\bydef \mathrm{Rep}(Q^{l},(1,\ldots,1))$, the representations of the unframed quiver, formed by deleting the vertex $\infty$. The corresponding quantised moment map is
\[
\hat\tau'(I^{(j)}) =X^{(j)}\partial^{(j)}  - X^{(j-1)}\partial^{(j-1)},
\] where $j=0,\ldots,l-1$. Note that this is $(-)^\mathrm{op}$-invariant. Let
\[
A_\chi \bydef \left(\frac{D_{V'}}{D_{V'}\langle \hat\tau'(A)-\chi(A)\, |\, A\in \mathfrak g\rangle}\right)^G.
\] Now $B_0,\ldots, B_{l-1},C,D$ generate $A_\chi$ and the relations are
\begin{align*} 
B_0\cdots B_{l-1} &= CD, &
 B_0-B_{l-1} &= \chi_0,&
 B_1-B_0 &=\chi_1,&
 &\ldots, &
 B_{l-1}-B_{l-2}&=\chi_{l-1}.
 \end{align*} Therefore, mapping $\chi_0$ to $-\sum_{i=1}^{l-1}\chi_i$ gives an isomorphism $A_\chi \cong D_{\hat\tau,\chi}^G$. 
 
 Now consider the results of \cite[Section 6.5]{BK}. They choose a basis $\langle u_1,\ldots,u_l\rangle =V$ so that $(\lambda_1,\ldots,\lambda_l)\in G$ acts on $(u_1,\ldots,u_l)$ by
\[
(\lambda_1\lambda_l^{-1}u_1,\ldots,\lambda_l\lambda_{l-1}^{-1}u_l).
\] They also choose the basis $\langle v_1,\ldots,v_{l}\rangle = \mathbb X(G)$, so that an arbitrary character is written $\phi = \sum_{i=1}^l \phi_i v_i$ for some $\phi_1,\ldots,\phi_l$. They then factor out by the diagonal action of $\CC^* $, let this group be denoted $\hat G$. The sublattice of characters such that $\sum_{i=1}^l \phi_i=0$ gives a basis for $\mathbb X(\hat G)$. They choose a new basis $w_i\bydef v_i-v_{i+1}$ for $i=1,\ldots,l-1$ so that a general character is written $\phi = \sum_{i=1}^l \hat{\chi}_i w_i$ where $\hat{\chi}_i \bydef \sum_{j=1}^i \phi_i$. The result is that the parameters $\hat{\chi}_1,\ldots,\hat{\chi}_{l-1}$ are related to the hyperplane parameters by the formula
\[
\hat{\chi}_i = h_i- h_0 + \frac{i-l}{l}
\] as $i=1,\ldots,l-1$.

This is converted into the notation used in this paper as follows. Identifying arbitrary elements $(t_0,\ldots,t_{l-1})$ and $(\lambda_1,\ldots,\lambda_l)$ by $t_i \leftrightarrow \lambda_{i+1}$ gives a $G$-equivariant
map between the vector spaces denoted $V$ by mapping $X_i\leftrightarrow u_{i+2}$, where the subscripts are taken modulo $l$. An arbitrary character, $\chi=\sum_{i=0}^{l-1} \chi_i \mathrm{tr}^{(i)}$, 
is identified with $\phi$ by $\chi_i = \phi_{i+1}$ for $i=0,\ldots,l-1$ and this gives 
$\chi_i = \hat{\chi}_{i+1}-\hat{\chi}_i$, for $i=0,\ldots,l-2$. 
Together, $\chi_i = h_{i+1} - h_{i} + \frac{1}{l} - \delta_{i,0}$, for $i=0,\ldots,l-2$, and, under the isomorphism between $A_\chi$ and $D_{\hat\tau,\chi}^G$, 
this gives $\chi_{l-1} = -\sum_{i=0}^{l-2} \chi_i = h_0 - h_{l-1} + \fr l$.
Converting the $h$'s to $k$'s by the formula $h_{i+1} - h_i = k_{1-i} - k_{-i}$, for all $i=0,\ldots,l-1$, gives  $\chi_i = k_{1-i}-k_{-i} + \fr l$, for $i=1,\ldots,l-1$.

Now, in order to convert the hyperplane paramaters to reflection parameters one uses the formula $k_{i+1}-k_i = \frac{-2}{l} \sum_{t=1}^{l-1} \zeta^{-it}c_t$, for all $i=1,\ldots,l-1$. The result is
 \[
 \chi_i = \frac{1}{l} \left( 1 -2\sum_{t=1}^{l-1} \zeta^{it}c_t\right).
 \]
 
 Now let $l>1$ and $n>1$, so that the results of \cite{GordonRemark} apply. Let $\tau$ be the quantised moment map chosen in that paper. As a consequence of having chosen $\tau$ to agree with a paper of Oblomkov, in that paper it is defined differently: as the differential of the $G$-action on $V$---the negative of that in \cite{GGS} which is defined as the differential of the action of $G$ on $\CC[V]$. It follows that $-\hat\tau$ and $\tau$ differ by a character of $\mathfrak{g}$. 
 \[ -\tau(I^{(i)}) =\begin{cases} -\sum_{s,t}X_{s,t}^{(i-1)}\partial_{s,t}^{(i-1)}+\sum_{s,t}X_{s,t}^{(i)}\partial_{s,t}^{(i)} &\textrm{ when $i\neq 0$,}\\
 -\sum_{s,t}X_{s,t}^{(-1)}\partial_{s,t}^{(-1)}+\sum_{s,t}X_{s,t}^{(0)}\partial_{s,t}^{(0)} - \sum_j x_j\partial_j &\textrm{ when $i=0$}.\end{cases}
\] Comparing this with $\hat\tau$ gives
\[
\hat\tau(I^{(i)})=\begin{cases} -\tau(I^{(i)}) & \textrm{ when $i\neq 0$,}\\ -\tau(I^{(0)})-\fr2 \mathrm{tr}(I_n) &\textrm{ when $i=0$;}
\end{cases}
\] so that $\hat\tau = -\tau - \fr2 \mathrm{tr}^{(0)}$. Comparing $D^G_{\hat\tau,\chi}$ with $D_{\tau,\chi_{c,k}}^G$ gives
\begin{align*}
D^G_{\hat\tau,\chi} \cong \left(\frac{D_V}{D_V\langle -\tau(A)-(\chi+\fr2 \mathrm{tr}^{(0)})(A)\, |\, A\in \mathfrak g\rangle}\right)^G \cong D^G_{\tau,-\left(\chi+\fr2\mathrm{tr}^{(0)}\right)}.
\end{align*} 
In \cite[Theorem 3.13]{GordonRemark}, he proves that $D^C_{\tau,\chi_{k,c}}\cong U_\cc$, where the character, $\chi_{k,c}$, is defined by 
\[
\chi_{k,c} \bydef (k+C_0)\mathrm{tr}^{(0)} + \sum_{i=1}^{l-1} C_i\mathrm{tr}^{(i)},
\] and the $C_i$'s and $k$ are related to $\cc$ by $k=-c_0$ and $C_i = \frac{1}{l}\left(1+2\sum_{t=1}^{l-1}\zeta^{it}c_t\right)-\delta_{i,0}$ for $i=0,\ldots,l-1$.

It follows that $\chi_0 = \fr2 - k - C_0$ and $\chi_i = -C_i$, for all $i=1,\ldots,l-1$. Together these give Equation \ref{G::equa:conversion}.
  \end{proof}
 \end{theo} Under the correspondence above, the parameter space for the spherical Cherednik algebra can now be thought of as $(\mathfrak g^* )^G$. Let $U_\chi$ denote the corresponding spherical rational Cherednik algebra.
 
The functor 
  \[
  \mathbb H\bydef \Gamma(X,-)^{F,G}\colon \left(\widetilde{\mathcal W}_X, F,G\right)^\mathsf{good}_\chi \longrightarrow {}_{U_\chi}\mathbf{mod}
  \] that takes $F$ and $G$-invariant global sections is called \textbf{quantum hamiltonian reduction}.
 Because $X$ is affine and $F$ and $G$ are reductive, $\mathbb H$ is exact.
 
\subsection{Constructing a W-Algebra on the GIT Quotient $Y_\theta$}
 
Recall that a choice of GIT parameter, $\theta\in \mathbb X(G)$, produces an open set of semistable points in $X$, denoted $X_\theta^\mathrm{ss}$. Assume that $\theta$ is chosen so that it does not lie on a GIT wall.
 The restriction functor $\mathsf{Res}\colon \mathcal W_X\longrightarrow \mathcal W_{X_\theta^\mathrm{ss}}$ is exact and induces an exact functor
 \[
 \mathsf{Res}\colon \left(\widetilde{\mathcal W}_X,F,G\right)^\mathsf{good}_\chi \longrightarrow \left(\widetilde{\mathcal W}_{X_\theta^\mathrm{ss}},F,G\right)^\mathsf{good}_\chi.
 \] 

The embedding functor, $I_\chi\colon \left(\mathcal W_X,G\right)_\chi\longrightarrow \left(\mathcal W_X,G\right)$, has a left adjoint, denoted 
 $\Phi_\chi\colon\left(\mathcal W_X,G\right)\longrightarrow \left(\mathcal W_X,G\right)_\chi$, defined by
 \[
 \Phi_\chi(\mathcal M)\bydef \frac{\mathcal M}{\left\langle \left\{A\cdot(-)-\hat\tau(A)+\chi(A)\, \middle|\, A\in \mathfrak g \right\}\right\rangle\mathcal M}.
 \]  Let $\mathcal L_{\chi}\bydef \Phi_\chi(\mathcal W_X)$. In \cite{KR}, they show that $\mathcal L_{\chi}$ is a good quasi-$G$-equivariant $\mathcal W_X$-module and that $\mathcal L_\chi$ is supported on the closed subset $\mu^{-1}(0)\subset X$.
 
Let $p\colon \mu^{-1}(0)\cap X_\theta^{\mathrm{ss}}\longrightarrow Y_\theta$ denote the GIT quotient map. Following \cite{KR}, define a sheaf of $\mathbf k$-algebras on $Y_\theta$ by,
\begin{equation}\label{W algebra on Y}
\mathcal W_{Y_\theta} \bydef \left( p_*  \mathcal{E}\mathit{nd}_{\mathcal W_{X^\mathrm{ss}}}(\mathcal L_\chi)^G \right )^\mathrm{op}.
\end{equation} In \cite[Proposition 2.8]{KR}, Kashiwara and Rouquier prove that $\mathcal W_{Y_\theta}$ is a W-algebra on ${Y_\theta}$, an $F$-action on $\mathcal W_{X^\mathrm{ss}}$ induces an $F$-action on $\mathcal W_{Y_\theta}$ and there is an equivalence of categories
\[
\mathbb E\colon\left(\widetilde{\mathcal W}_{X^\mathrm{ss}},F,G\right)^\mathsf{good}_\chi \stackrel{\simeq}{\longrightarrow} \left(\widetilde{\mathcal W}_{Y_\theta},F\right)^\mathsf{good}.
\]

By Proposition \ref{G::prop:codim at least 2} the semistable points have codimension at least two. Hartogs' Extension Theorem gives $\Gamma(X,\widetilde{\mathcal W}_X)^{F,G}\cong \Gamma(X^\mathrm{ss},\widetilde{\mathcal W}_{X^\mathrm{ss}})^{F,G}$.
 Because $\mathbb E$ is an equivalence, taking $F$-invariants of global sections gives a functor
\[
\mathbb T\bydef \Gamma(X,-)^F\colon \left(\widetilde{\mathcal W}_{Y_\theta},F\right)^\mathsf{good}\longrightarrow {}_{U_\chi}\mathbf{mod}.
\] Together, these functors fit into the following (not necessarily commutative) diagram.
 \begin{equation}\label{G::diag:localisation}
 \parbox{1ex}{\xymatrix@C=10ex{
 (\widetilde{\mathcal W}_{X}, F, G)^\mathsf{good}_\chi \ar[d]^{\mathsf{Res}} \ar@(r,u)[ddr]^{\mathbb H}& \\ 
 (\widetilde{\mathcal W}_{X^\mathrm{ss}_\theta}, F, G)^\mathsf{good}_\chi \ar[d]^{\mathbb E} & \\
 (\widetilde{\mathcal W}_{Y_\theta}, F)^\mathsf{good} \ar[r]^-{\mathbb T} & {}_{U_\chi}\mathbf{mod}
 }}
 \end{equation}

 In \cite{BPW}, when $\mathbb T$ induces a derived equivalence they say that \textbf{derived localisation holds for $\chi$}. When $\mathbb T$ is an 
equivalence of abelian categories they say that \textbf{localisation holds for $\chi$}. 

 \begin{defi} \label{G::defi:good.n.bad.params} Say that a character, $\chi\in (\mathfrak g^* )^G$, is \textbf{bad for $\theta\in \mathbb X(G)$} if
 $\ker\mathsf{Res}\nsubseteq \ker\mathbb H$. Say that $\chi$ is \textbf{bad} if it is bad for all $\theta$ and \textbf{good} if it is not bad;
 that is, there exists some $\theta$ such that $\ker \mathsf{Res}\subseteq \ker\mathbb H$.
 \end{defi}
If a character $\chi$ is good for some $\theta$ and the corresponding parameter $\cc$ is spherical then there are several nice consequences. McGerty and Nevins, in \cite[Corollary 7.5]{MN} and \cite[Lemma 3.9]{MN}, prove that when $\chi$ is spherical, derived localisation holds for $\chi$. This result is used to prove the following theorem.

 \begin{theo}\label{G::theo:Commutes} Suppose that $\chi\in (\mathfrak g^* )^G$ is good for $\theta\in \mathbb X(G)$ and that the corresponding value of $\cc$ given by Theorem \ref{G::theo:chi and params} is spherical. Then 
\begin{enumerate}[(i)]
\item Diagram \ref{G::diag:localisation} commutes,
\item localisation holds for $\chi$,
\item the kernel of $\mathbb H$ is precisely those sheaves supported on $\us$, that is,
\[
\ker {\mathbb H} = \ker \mathsf{Res}.
\]
\end{enumerate}
\begin{proof}

First, if $\ker \mathsf{Res}\subseteq \ker \mathbb H$, then $\mathbb H=\Gamma(X,-)^{F,G}$ factors through $\mathsf{Res}$ as 
\[
\mathbb H=\Gamma(X^\mathrm{ss},-)^{F,G}\circ \mathsf{Res}=\mathbb T \circ \mathbb E\circ \mathsf{Res}.
\]

Secondly, I claim that if Diagram \ref{G::diag:localisation} commutes then the functor $\mathbb T$ is exact. Since $\mathbb T$ is the $F$-invariant global sections functor on the space $Y_\theta$, it is automatically left exact; so it is sufficient to show that it is right exact. 

The restriction functor, $\mathsf{Res}$, has a left adjoint, $\mathsf{Res}_!$, such that $\mathsf{Res}\circ \mathsf{Res}_! \cong \mathrm{id}$. Indeed, following \cite{BPW}, define the \textbf{Kirwan functor},
\[
\kappa \colon \left( \widetilde{\mathcal W}_X,F,G\right)^\mathsf{good}\longrightarrow \left( \widetilde{\mathcal W}_{Y_\theta},F\right)^\mathsf{good},
\] by $\kappa(\mathcal N)\bydef p_* \mathcal{H}\mathit{om}\left(\mathcal L_\chi, \mathcal N|_{X^\mathrm{ss}}\right)$.
In \cite[Lemma 5.18]{BPW}, they show that it has a left adjoint $\kappa_!$ such that $\kappa\circ \kappa_!\cong \mathrm{id}$. 

Recall that the forgetful functor $I_\chi\colon (\widetilde{\mathcal W},F,G)^\mathsf{good}_\chi\longrightarrow (\widetilde{\mathcal W},F,G)^\mathsf{good}$ also has a left adjoint: 
$\Phi_\chi$ and it satisfies $I_\chi \circ \Phi_\chi\cong \mathrm{id}$. Now define
\[
\mathsf{Res}_! \bydef \Phi_\chi\circ \kappa_! \circ \mathbb E.
\] Being the composition of two left adjoints and an equivalence, it is a left adjoint and, because Diagram \ref{G::diag:localisation} commutes,
\[
\mathsf{Res}\circ\mathsf{Res}_! = \mathsf{Res}\circ \Phi_\chi\circ \kappa_! \circ \mathbb E \cong \mathbb E^{-1} \circ \kappa \circ \kappa_! \circ \mathbb E \cong \mathrm{id}.
\] Now, $\mathbb T\circ \mathbb E$ is right exact because
\[
\mathbb T\circ \mathbb E \cong \mathbb T\circ \mathbb E \circ \mathsf{Res}\circ \mathsf{Res}_! = \mathbb H \circ \mathsf{Res}_!,
\] the composition of two right exact functors. It follows that $\mathbb T\circ \mathbb E$, and therefore $\mathbb T$, is right exact. This completes the proof of the claim.

The exactness of $\mathbb T$ now implies that $\mathbb T$ is an equivalence. Indeed, let $\mathbb S\bydef \widetilde{\mathcal W}_{Y_\theta} \otimes_{U_\chi} -$, the left adjoint of $\mathbb T=\Gamma(Y_\theta,-)^F$. 
Being a left adjoint, $\mathbb S$ is right exact. Let $R\mathbb T$ denote the right derived functor of $\mathbb T$ and $L\mathbb S$ the left derived functor of $\mathbb S$. 
Because $\mathbb T$ is exact, as derived functors, $R\mathbb T=\mathbb T$. Let $\mathcal M\in (\widetilde{\mathcal W}_{Y_\theta}, F)^\mathsf{good}$. 

By \cite[Corollary 7.5]{MN}, when $\chi$ is spherical, $R\mathbb T$ and $L\mathbb S$ are equivalences so
\[
L\mathbb S\circ R\mathbb T\mathcal M=\widetilde{\mathcal W}_{Y_\theta}\stackrel{L}{\otimes}_{U_\chi}\Gamma(Y_\theta,\mathcal M)^F\longrightarrow \mathcal M
\] is a quasi-isomorphism. Therefore,
\[
\mathbb S\circ \mathbb T\mathcal M=\widetilde{\mathcal W}_{Y_\theta}\otimes_{U_\chi}\Gamma(Y_\theta,\mathcal M)^F\longrightarrow \mathcal M
\] is surjective.

Now, $\mathbb T\mathcal M$ is finitely generated, so choosing $m$ generators gives a surjective map $U_\chi^{\oplus m}\longrightarrow \mathbb T\mathcal M$. Applying the right exact functor $\mathbb S$ gives a surjection which, when composed with the surjection above, gives a surjection
\[
\mathbb S U_\chi^{\oplus m}\cong \widetilde{\mathcal W}_{Y_\theta}^{\oplus m} \longrightarrow \mathbb S\circ \mathbb T \mathcal M \longrightarrow \mathcal M.
\]

This shows that any module is a quotient of some power of $\widetilde{\mathcal W}_{Y_\theta}$, in particular, it can be applied to the kernel of the composition of the maps above to give an exact sequence,
\[
\xymatrix{\widetilde{\mathcal W}_{Y_\theta}^{\oplus m'} \ar[d] \ar[r] & \widetilde{\mathcal W}_{Y_\theta}^{\oplus m} \ar@{->>}[r] \ar[d]& \mathcal M \ar[d] \\
\mathbb S\circ \mathbb T \widetilde{\mathcal W}_{Y_\theta}^{\oplus m'} \ar[r]& \mathbb S\circ \mathbb T\widetilde{\mathcal W}_{Y_\theta}^{\oplus m} \ar@{->>}[r] & \mathbb S\circ \mathbb T\mathcal M,}
\] where the vertical maps are from the natural transformation between $\mathbb S\circ \mathbb T$ and the identity functor. The first two of these are isomorphisms induced from
\[
\widetilde{\mathcal W}_{Y_\theta} \otimes_{U_\chi}U_\chi \cong \widetilde{\mathcal W}_{Y_\theta},
\] which implies the third vertical map is an isomorphism.
Now $\mathbb S$ and $\mathbb T$ are inverse equivalences, so localisation holds for $\chi$.

Finally, since $\mathbb T \circ \mathbb E=\Gamma(X^\mathrm{ss},-)^{F,G}$ is now an equivalence, $\ker \mathsf{Res}=\ker (\mathbb T\circ \mathbb E \circ \mathsf{Res}) = \ker \mathbb H$.
\end{proof}
\end{theo}

\section{The Kirwan--Ness Stratification and the McGerty--Nevins Criterion}

Let $T\cong(\CC^* )^n$ be a maximal torus in $G$ and $\mathbb Y(T)$ the group of one-parameter subgroups of $T$. Let $\mathbb{Y}(T)_{\R}$ denote the tensor product $\mathbb{Y}(T)\otimes_{\Z} \R$. 
The Weyl group, $W = N_G(T)/T$, acts on $\mathbb{Y}(T)$ and $\mathbb{Y}(T)_{\R}$ by conjugation. 
Choose a $W$-invariant inner product $(-,-)$ on $\mathbb{Y}(T)_{\R}$ and denote the 
associated norm by $|| \cdot ||_T$.
When the image of $\lambda\in \mathbb Y(G)$ lies inside $T$ one can use the isomorphism 
\[
\mathbb Y(T)\stackrel{\cong}{\longrightarrow} \Z^n;\qquad \lambda\mapsto d\lambda(1)
\]
to write $\lambda$ as an $n$-tuple of integers. In order to simplify notation, the symbol $\lambda$ will simultaneously stand for the map $\CC^*  \rightarrow G$ and the $n$-tuple of integers $d\lambda(1)$. 

Let $x\in X$ and suppose there exists a $\lambda\in \mathbb Y(G)$ such that $\lim_{t\rightarrow 0}\lambda(t)\cdot x = x_0\in X$ exists. Then $\lambda$ induces a $\CC^* $-action on the fibre, $\mathcal L_{x_0}$, given by $t\cdot(x_0,l)=(x_0,t^rl)$ for $t\in \CC^* $, $l\in \mathcal L_{x_0}$ and some $r\in \Z$. Let 
$\mu(x,\lambda)$ denote this integer, $r$. If the limit does not exist then define $\mu(x,\lambda)=\infty$. The \textbf{Hilbert--Mumford criterion} states that a point, $x\in X$, is unstable if and only if there exists a one-parameter subgroup, $\lambda$, such that $\mu^\theta(x,\lambda)<0$. See \cite[Proposition 2.5]{Hosk}. Define
\[
X_\lambda^\theta \bydef \{ x\in X\, |\, \mu^\theta(x,\lambda)<0\}.
\]
One can also use the Hilbert-Mumford criterion to measure the extent to which a point is unstable. 
Let $\lambda\in \mathbb{Y}(G)$ and pick a $g$ such that $\mathrm{Ad}(g)\cdot \lambda \in \mathbb{Y}(T)$. Define a norm of $\lambda$ by
\[
\|\lambda \| \bydef \|\mathrm{Ad}(g)\cdot \lambda\|_T
\]  
and a function
\[
M^\theta(x) \bydef \sup_{\lambda\in\mathbb{Y}(G)\setminus \{0\}} \left\{\frac{-\mu(x,\lambda)}{\|\lambda\|}\right\}.
\]

\begin{defi}
A one-parameter subgroup is called \textbf{primitive} (or \textbf{indivisible}) if it cannot be written as a positive multiple of another.
A primitive one-parameter subgroup, $\lambda$, is called $\theta$-\textbf{optimal} (or just \textbf{optimal}) for 
$x$ if $M^\theta(x)$ realises its supremum, $\frac{-\mu^\theta(x,\lambda)}{\|\lambda\|}$, at $\lambda$.
Finally, let $\lambda$ and $\mu$ be one-parameter subgroups which are not a positive integer multiple of one another. Say that $\lambda$ \textbf{$\theta$-dominates} $\mu$ 
if $X^\theta_\mu \subseteq X^\theta_\lambda$ and $\frac{\mu^\theta(x,\mu)}{\|\mu\|}\ge \frac{\mu^\theta(x,\lambda)}{\|\lambda\|}$ for all $x\in X_\mu^\theta$.
\end{defi} It follows that optimal subgroups are not 
dominated by any other. Let $\Gamma_\theta$ denote the set of all optimal one-parameter subgroups for the stability condition $\theta$. 
The set $\Gamma_\theta$ is complicated: it varies even as one varies $\theta$ inside a single GIT chamber. 
The following theorem was proved independently by Kirwan and Ness; see \cite[Sections 12-13]{Kirwan}
and \cite{Ness} respectively.

\begin{theo}(Kirwan--Ness) The unstable locus, $\us$, has the following \textbf{Kirwan--Ness stratification}. 
\begin{itemize}
\item[(KN1)] There is a decomposition into non-empty smooth locally closed strata 
\[X^{\mathrm us} = \coprod_{d, \langle \lambda \rangle} {S}^+_{d, \langle \lambda \rangle},\] 
where $d$ is a positive real number and $\langle \lambda \rangle$ a conjugacy class of one-parameter subgroups of $G$.
\item[(KN2)] There is an enumeration of the one-parameter subgroups appearing in (KN1) 
by representatives $\lambda_1, \ldots, \lambda_p\in \mathbb{Y}(T)$, such that 
$\overline{S}^+_{d,\langle \lambda_i\rangle} \cap S^+_{d',\langle \lambda_j \rangle}\neq \emptyset$ only if $i<j$ and $d<d'$.
\item[(KN3)] For $1\leq i \leq p$, set 
\[
P(\lambda_i) \bydef \left\{ g\in G\, \middle|\, \lim_{t\to 0} \lambda_i(t)g\lambda_i(t)^{-1} \textrm{ exists }\right\}.
\] 
There exists a smooth, locally-closed $P(\lambda_i)$-stable subvariety, $S_{d, \lambda_i}$, of $X$, such that the action map 
induces 
\[
G \stackrel{P(\lambda_i)}{\times} S_{d,\lambda_i} \cong S^+_{d,\langle \lambda_i \rangle}.
\]
\item[(KN4)] For $1\leq i \leq p$, set 
\[
Z_{d,\lambda_i} \bydef \left\{ x\in S_{d,\lambda_i}\, \middle|\, \lambda_i(t) \cdot x = x \textrm{ for all }t\in \mathbb G_m \right\},
\] 
and $Z_G(\lambda_i) \bydef \left\{g\in G\, \middle|\,  \lambda_i(t)g\lambda_i(t)^{-1} = g\textrm{ for all }t\in \mathbb G_m\right\}$. Then the variety, $Z_{d,\lambda_i}$, 
is a $Z_G(\lambda_i)$-stable, smooth, locally closed, subvariety of $X$ such that 
\[
S_{d,\lambda_i} = \left\{ x\in X\, \middle|\, \lim_{t\to 0} \lambda_i(t) \cdot x \in Z_{d,\lambda_i}\right\}.
\]
\item[(KN5)] For $1\leq i \leq p$, let 
\[
Z_{d,\lambda_i} \bydef \coprod Z_{d,\lambda_i, j}
\] 
be the decomposition of $Z_{d,\lambda_i}$ into connected components. For each $j$ that appears, the morphism,  
\[
p_{d,\lambda_i, j}\colon S_{d,\lambda_i, j} \longrightarrow Z_{d, \lambda_i, j}; \quad x \mapsto
\lim_{t\to 0} \lambda_i(t)\cdot x,
\] 
is a locally trivial fibration by affine spaces. 
\end{itemize} 
\end{theo}

The strata can be described in terms of unstability by
\[
S^+_{d,\langle \lambda \rangle} =\left\{x\in X\, \middle|\, M(x) = d \textrm{ and there exists a } g\in G 
\textrm{ with }\mathrm{Ad} g\cdot \lambda \textrm{ optimal for } x\right\}
\] 
and 
\[
S_{d,\lambda} = \left\{x\in X\, \middle|\, M(x) = d \textrm{ and } \lambda \textrm{ is optimal for }x\right\}.
\]

\begin{rema} The index, $d$, for each Kirwan--Ness stratum is redundant in the rest of the paper because, when $X$ is affine, it is determined by $\lambda$ and $\theta$. That is,
\[
d=-\frac{\theta\cdot \lambda}{\|\lambda\|}.
\] The subscript $d$ will be dropped from the notation from now on.
\end{rema}

\subsection{On the Kirwan--Ness Stratification for $T^* \mathrm{Rep}(Q^l_\infty,\epsilon)$} \label{G::sect:gl1n.strata}

Let $\lambda= \left(a_k^{(m)}\right)_{k,m}\in \mathbb Y(T)$ and $\theta=(\theta^0,\ldots,\theta^{l-1})$, so that for some $x\in X$ such that $\lim_{t\rightarrow 0}
 \lambda(t)\cdot x$ exists,
 \[
 \mu^\theta(x,\lambda) = \sum_{m,k}\theta^ma_k^{(m)}.
 \]
 When it is understood from context which stability condition is being used, $\theta$ will be dropped from all the above notation. 
 
\begin{lemm}\label{G::lemm:reording improves} Fix positive integers, $x_1\le\cdots\le x_n\in \N$ and $y_1\le\cdots\le y_n\in \N$, for some $n\ge 2$. The function 
  \[
  f\colon \mathfrak S_n\longrightarrow \N;\quad \sigma \longmapsto \sum_{i=1}^n x_{\sigma(i)}y_i
  \] is maximal at $\mathrm{id}\in \mathfrak S_n$. 
 \begin{proof} This is obvious for $n=2$. An induction then proves the lemma.
 \end{proof}
 \end{lemm}
 
 \begin{prop}\label{G::prop:codim at least 2} Let $X$, $G$ and $\theta$ be as above. The locus of unstable points, $\us$, have codimension at least two in $X$.
\begin{proof}
Since there are only finitely many strata, there exists some optimal one-parameter subgroup, $\lambda$, such that $\dim S_{\langle \lambda \rangle}^+ = \dim \us$. By (KN3) each strata can be decomposed as
\[
S_{\langle \lambda \rangle}^+ \cong G \stackrel{P(\lambda)}{\times} S_\lambda.
\] By counting dimensions, $\mathrm{codim}_X\us = \dim X- \dim S_\lambda - \dim G + \dim P(\lambda)$.

First, consider the case, $n=1$. Here, $\dim G - \dim P(\lambda)=0$, so the codimension is the number of negative weights of $\lambda$ acting on $T^*V$. This is the space of representations of the unframed quiver, $\mathrm{Rep}(\overline Q^l,\gamma)$ with the action of 
$\hat G= (\CC^* )^l/\CC^* $. As an element of $\mathbb Y(G)$, an optimal one-parameter subgroup is an $l$-tuple, $\lambda=(a_0,\ldots,a_{l-1})$. It must be non-zero in $\mathbb Y(\hat G)$, so there is some $i\neq j$ such that $a_i\neq a_{i-1}$ and $a_j\neq a_{j-1}$. This means that $X^{(i)}$ and $X^{(j)}$ are non-fixed eigenvectors in $V$, which implies that there are at least two negative weights of $\lambda$ acting on $T^* V$, this implies that $X^\mathrm{ss}$ contains a subspace of dimension two.

Now let $n\ge 2$. Let $\mathbf X^\lambda$ denote the subspace of $\mathbf X=\mathbf X^{(0)}\oplus \cdots \oplus \mathbf X^{(l-1)}$ fixed by $\lambda$. Let $G^\lambda$ denote the subspace of $G$ fixed by $\lambda$ (acting by conjugation). Note that
\[
2\left(\dim G - \dim P(\lambda)\right) =  ln^2 - \dim (G^\lambda).
\]
I claim that $\dim (\mathbf X^\lambda)\le \dim (G^\lambda)$.

Let $\lambda\bydef (\lambda^{(0)},\ldots,\lambda^{(l-1)})$, where, for each $k=0,\ldots,l-1$, $\lambda^{(k)}\bydef (a_1^{(k)},\ldots,a_n^{(k)})\in \Z^{n}$. By conjugating this by the action of the Weyl group, $\mathfrak S_n \times \cdots \times \mathfrak S_n$, assume that the entries within each component are in ascending order; that is, for all $k=0,\ldots,l-1$ and $i=1,\ldots,n-1$, $a_i^{(k)}\le a_{i+1}^{(k)}$. Rename these entries so that, for each $k=0,\ldots,l-1$, after removing duplicate entries, $\lambda^{(k)}$ would look like, $(b_1^{(k)},\ldots,b_{p_k}^{(k)})$ with $b_i^{(k)}<b_{i+1}^{(k)}$ for all $i=1,\ldots,p_k-1$. Let $n_i^{(k)}$ be the number of times $b_i^{(k)}$ appears in $\lambda^{(k)}$. With this new notation each component of the subgroup looks like
\[
\lambda^{(k)}=(\underbrace{b_1^{(k)},\ldots,b_1^{(k)}}_{n_1^{(k)}},\underbrace{b_2^{(k)},\ldots,b_2^{(k)}}_{n_2^{(k)}},\ldots, \underbrace{b_{p_k}^{(k)},\ldots,b_{p_k}^{(k)}}_{n_{p_k}^{(k)}}).
\] Note that, for each $k=0,\ldots,l-1$, $n_1^{(k)}+\cdots+n_{p_k}^{(k)}=n$ partitions $n$. Now, for $k=0,\ldots,l-1$, the weights of $\lambda^{(k)}$ acting on $\mathrm{GL}_{n_1^{(k)}+\cdots+n_{p_k}^{(k)}}(\CC)$ are 
$\begin{pmatrix} a_i^{(k)}-a_j^{(k)} \end{pmatrix}_{i,j}$, and these are zero precisely on the square blocks, cut out by the $n_i^{(k)}$, that run down the diagonal. That is to say, $\dim (G^\lambda)=\sum_{i=1}^{p_k}\left(n_i^{(k)}\right)^2$.

On the other hand, $\lambda$ acts on $\mathbf X^{(k)}$ with weights $A\bydef \begin{pmatrix} a^{(k+1)}_i - a^{(k)}_j \end{pmatrix}_{i,j}$.
Partition the rows of this matrix into $n=n_1^{(k+1)}+\cdots+n_{p_{k+1}}^{(k+1)}$ and the columns into $n=n_1^{(k)}+\cdots+n_{p_k}^{(k)}$. This divides the matrix of $\lambda$-weights into rectangular blocks inside each of which the weight is constant. Form a new $p_{k+1}\times p_k$ matrix, $B$, by treating each rectangular block of $A$ as a single entry: $B\bydef\begin{pmatrix} b_i^{(k+1)}-b_j^{(k)}\end{pmatrix}_{i,j}$. Because $\lambda^{(k)}$ and $\lambda^{(k+1)}$ are assumed to be increasing, each column of $B$ is strictly increasing as one moves down the column and each row is strictly decreasing as one moves left-to-right along the row. Clearly then, each row or column can contain at most one zero. For each $(i,j)$ such that $B_{i,j}=0$, $A$ contains exactly $n_j^{(k)}n_i^{(k+1)}$ zeroes.

For each $k=0,\ldots,l-1$, put the positive integers, $\left\{n_1^{(k)},\ldots,n_{p_k}^{(k)}\right\}$, into increasing order and rename them $\{m_1^{(k)},\ldots,m_{p_k}^{(k)}\}$.
Suppose, without loss of generality that $p_{k}\ge p_{k+1}$. Applying Lemma \ref{G::lemm:reording improves} to the two lists $\{m^{(k)}_1,\ldots,m_{p_k}^{(k)}\}$ and $\{m^{(k+1)}_1,\ldots,m_{p_{k+1}}^{(k+1)}\}$ shows that the number of zeroes of $\lambda$ acting on $\mathbf X^{(k)}$ is bounded above by
\[
 \sum_{i=p_k-p_{k+1}+1}^{p_k} m_i^{(k)}m_{i-p_k+p_{k+1}}^{(k+1)}.
\] When $p_k\le p_{k+1}$ the dummy variable runs from $p_{k+1}-p_k+1$ to $p_{k+1}$ and the subscript of $m^{(k)}$ and $m^{(k+1)}$ are adjusted appropriately. It follows that
\[
\dim(\mathbf X^\lambda) \le \sum_{k\, |\, p_k\ge p_{k+1}}^{l-1}\sum_{i=p_k-p_{k+1}+1}^{p_k} m_i^{(k)}m_{i-p_k+p_{k+1}}^{(k+1)}+ \sum_{k\, |\, p_k< p_{k+1}}^{l-1}\sum_{i=p_{k+1}-p_{k}+1}^{p_{k+1}} m_{i-p_{k+1}+p_{k}}^{(k)}m_{i}^{(k+1)}.
\] But now note that $\{m_i^{(k)}\, |\, k=0,\ldots,l-1,\, i=1,\ldots,p_k\}=\{n_i^{(k)}\, |\, k=0,\ldots,l-1,\, i=1,\ldots,p_k\}$. Applying Lemma \ref{G::lemm:reording improves} again to two copies of this larger set gives
\[
\sum_{k=0}^{l-1}\sum_i m_i^{(k)}m_i^{(k+1)} \le \sum_{k,i} \left (n_i^{(k)}\right)^2 = \dim (G^\lambda).
\] This proves the claim.

Let $\dim (\mathbf X\oplus \mathbf Y)_-$ and $\dim (\mathbf v\oplus \mathbf w)_-$ denote the number of negative weights of $\lambda$ acting on $\mathbf X\oplus \mathbf Y$ and $\mathbf v\oplus \mathbf w$ respectively. Since the action of $\lambda$ is hamiltonian, the weights of $\lambda$ on $\mathbf X$ are the negatives of those on $\mathbf Y$. That means $\dim (\mathbf X\oplus \mathbf Y)_- + \dim (\mathbf X^\lambda) = ln^2$.

Together,
\begin{align*}
\mathrm{codim}_X\us &= \dim X- \dim S_{\langle \lambda \rangle}^+\\
			&= \dim X- \dim S_\lambda - \dim G + \dim P(\lambda)\\
			& = \dim (\mathbf X\oplus \mathbf Y)_-+\dim (\mathbf v\oplus \mathbf w)_-- \dim G + \dim P(\lambda)\\
			& = ln^2 - \dim (\mathbf X^\lambda) - \dim G + \dim P(\lambda) + \dim (\mathbf v\oplus \mathbf w)_-\\
			& \ge ln^2 - \dim (G^\lambda) - \dim G + \dim P(\lambda) + \dim (\mathbf v\oplus \mathbf w)_-\\
			& = \dim G - \dim P(\lambda) + \dim (\mathbf v\oplus \mathbf w)_-.
\end{align*} This is a sum of two non-negative integers.

Now consider the worst case scenarios. If $\dim (\mathbf v\oplus \mathbf w)_-=1$ then $\lambda^{(0)}=(0,\ldots,0,1)$, and the contribution from this component of $\lambda$ gives $\dim G-\dim P(\lambda)\ge 1$.
If $\dim (\mathbf v\oplus \mathbf w)_-=0$ then $\lambda^{(0)}=0$. If $\dim G - \dim P(\lambda)=0$ then $\mathrm{codim}_X\us=\dim X- \dim S_\lambda$ and the argument follows the case $n=1$ above. Otherwise, suppose $\dim G - \dim P(\lambda)=1$. There is only one possible form that $\lambda$ can take now. First, $\lambda$ must act trivially on all but one component of $G$, the $k^\textrm{th}$ say. The block of the matrix of weights of $\lambda^{(k)}$ acting on $\mathrm{GL}_n(\CC)$ that contributes $1$ to $\dim G - \dim P(\lambda)$ must be one-by-one, which implies, using the notation above that $n_1^{(k)}=n_2^{(k)}=1$. Therefore, $n=2$ and $a_1^{(k)}\neq a_2^{(k)}$, but for all $k'\neq k$, $a_1^{(k')}=a_2^{(k')}$. The weights of $\lambda$ acting on $\mathbf X^{(k-1)}\times \mathbf X^{(k)}$ are
\[
\begin{pmatrix} a_1^{(k)}-a_1^{(k-1)} & a_1^{(k)}-a_1^{(k-1)} \\ a_2^{(k)}-a_1^{(k-1)} & a_2^{(k)}-a_1^{(k-1)} \end{pmatrix} , \begin{pmatrix} a_1^{(k+1)}-a_1^{(k)} & a_1^{(k+1)}-a_2^{(k)} \\ a_1^{(k+1)}-a_1^{(k)} & a_1^{(k+1)}-a_2^{(k)} \end{pmatrix}.
\] At least four of these are non-zero, so there are at least four negative weights of $\lambda$ acting on $\mathbf X^{(k-1)}\times \mathbf X^{(k)} \times \mathbf Y^{(k-1)}\times \mathbf Y^{(k)}$. Therefore, $\dim X- \dim S_\lambda - \dim G + \dim P(\lambda)\ge 4-1=3$.
\end{proof}
\end{prop}

 \subsubsection{The Kirwan--Ness Stratification for $\theta=(-1,\ldots,-1)\in \Z^l$}
 
 By \cite[Lemma 4.3]{Gordon}, the GIT parameter $\theta=(-1,\ldots,-1)$ never lies on a GIT wall. 
 The aim of this subsection is to find the set, $\Gamma_\theta$, of optimal one-parameter subgroups for this stability condition.
 The adjoint action of $N(T)$ on $\mathbb Y(T)$ factors through the Weyl group $N(T)/T\cong (\mathfrak S_n)^l$. This permutes the entries of $\lambda=\left(a_k^{(m)}\right)_{k,m}\in \mathbb Y(T)$
 so that, for an $l$-tuple of permutations, $(\sigma^{(m)})\in (\mathfrak S_n)^l$, $(\sigma^{(m)})\cdot \lambda = \left(a_{{\sigma^{(m)}}^{-1}(k)}^{(m)}\right)$.
 Define
 \[
 I\bydef \left\{\left. \left(a_k^{(m)}\right)_{k,m}\, \right|\, a_k^{(m)}\in\{0,1\}\right\}\setminus\{(0,\ldots,0)\}.
 \]
 
 \begin{lemm}
  Let $\lambda=\left(a_k^{(m)}\right)_{k,m}\in \mathbb Y(T)$ and $\nu=\left(b_k^{(m)}\right)_{k,m}$ where $b_k^{(m)}\bydef\begin{cases} 0&\textrm{ if }a_k^{(m)}\le 0\\ 1&\textrm{ if }a_k^{(m)}>0\end{cases}$.
 Then $\nu$ dominates $\lambda$.
 \begin{proof}
 Suppose that there exists an $x\in X_\lambda\setminus X_\nu$. Then, either $\mu(x,\nu)=-\sum b_k^{(m)}\ge 0$, in which case $\mu(x,\lambda)=-\sum a_k^{(m)}\ge - \sum b_k^{(m)}= 0$, contradicting $x\in X_\lambda$, or $\lim_{t\rightarrow 0} \nu(t)\cdot x$ doesn't exist. 
 
 It follows that one of the coordinates of $x$ has a negative weight with respect to $\nu$. Suppose that it is $\mathbf X_{ij}^{(m)}$. Then $b_i^{(m+1)}< b_j^{(m)}$ so $b_i^{(m+1)}=0$ and $b_j^{(m)}=1$ which implies that
 $a_i^{(m+1)} \le 0 < a_j^{(m)}$ which contradicts $x$ being unstable for $\lambda$. A similar argument shows that 
 $x$ cannot be unstable at a coordinate $\mathbf Y^{(m)}_{ij}$ for $\lambda$ unless it is unstable for $\nu$.
 Note that $\mathbf v_i$ cannot have a negative $\nu$-weight.
 Suppose that $\mathbf w_i$ has a negative $\nu$-weight. Then $b_i^{(0)}=1$, so $a_i^{(0)}>0$ and $x$ is not unstable for $\lambda$. This contradiction implies
 $X_\lambda\subseteq X_\nu$.
 
 It remains to prove that $\frac{\mu(x,\nu)}{\|\nu\|} \le \frac{\mu(x,\lambda)}{\|\lambda\|}$. This is done by two claims. First, increasing the negative entries of $\lambda$ to zero 
 decreases the value of this ratio. Second, once all the entries of $\lambda$ are non-negative, changing those which are non-zero to one doesn't increase the value of the ratio.
 
 Suppose $a_k^{(m)}<0$ for some $k,m$. Let $\lambda'$ be the one-parameter subgroup whose entries agree with $\lambda$ everywhere except at the $(k,m)$ position where it is zero; then $\|\lambda'\|<\|\lambda\|$. Now, 
\[
 \frac{\mu(x,\lambda)}{\|\lambda\|}= \frac{\mu(x,\lambda')}{\|\lambda\|}-\frac{a_k^{(m)}}{\|\lambda\|}  > \frac{\mu(x,\lambda')}{\|\lambda\|} >\frac{\mu(x,\lambda')}{\|\lambda'\|},
 \] where the last inequality holds because $\mu(x,\lambda')$ is negative. By repeating this argument for each negative entry, the first claim is proved.
 
 To prove the second claim, suppose that $\lambda=\left(a_k^{(m)}\right)_{k,m}$ is such that $\lim_{t\rightarrow 0}\lambda'(t)\cdot x$ exists and $a_k^{(m)}\ge 0$ for all $k,m$. 
 Let $N$ be the number of non-zero coordinates. Then, by the Cauchy-Schwartz inequality,
 \begin{align*}
  \left(\sum_{k,m} a_k^{(m)}.1\right)^2 &\le \left( \sum_{k,m} {a_k^{(m)}}^2\right)\left(\sum_{i=1}^N 1^2\right).\\
 \intertext{Since all the coordinates are non-negative, taking square roots gives}
 \sum_{k,m} a_k^{(m)} &\le \sqrt{ \sum_{k,m} {a_k^{(m)}}^2}\sqrt{N};\\
 \intertext{so that}
 \frac{\mu(x,\lambda)}{\|\lambda\|}= \frac{-\sum a_k^{(m)}}{\sqrt{ \sum {a_k^{(m)}}^2}} &\le -\sqrt N = \frac{\mu(x,\nu)}{\|\nu\|}.
 \end{align*} 
 \end{proof}
 \end{lemm}

 \begin{coro}\label{coro:G.zerosandones}
  If $\lambda\in \mathbb Y(T)$ is optimal for $x\in X_\lambda$ then $\lambda\in I$.
 \begin{proof}
  It suffices to check that the $b_k^{(m)}$ constructed in the lemma could not all be zero; if this were the case then
 $a_k^{(m)} \le 0$ for all $k$ and $m$ so $\mu(\lambda,x)\nless 0$: a contradiction.
 \end{proof}
 \end{coro}
 
 Next, the set of candidate optimal subgroups is reduced by removing those which are dominated by another.
 Define a relation, $\rightarrow$, on the set of pairs $\{1,\ldots,n\}\times\Z_l$ by
 \[
  (i,m)\rightarrow \begin{cases} (j,m-1) &\textrm{ if and only if }\mathbf X_{ij}^{(m-1)}\neq 0 \\ (j,m+1)&\textrm{ if and only if }\mathbf Y_{ij}^{(m)}\neq 0 \end{cases}.
 \] Define a second relation, $\rightsquigarrow$, to be the transitive closure of $\rightarrow$; that is,
 \[
  (i,m)\rightsquigarrow (j,m') \quad \Longleftrightarrow \parbox{7cm}{\centering there exists a sequence of pairs, $((i_p,m_p)\, |\, p=1,\ldots,r)$, such that 
 $(i,m)\rightarrow (i_0,m_0)\rightarrow\cdots\rightarrow (i_r,m_r)\rightarrow (j,m')$.}
 \]
 
 \begin{lemm}\label{lemm:G.rightsquigarrow} Suppose that $x\in X_\lambda$ for some $\lambda = \left(a_k^{(m)}\right)_{k,m}$. If $(i,m)\rightsquigarrow (j,m')$ then $a_i^{(m)}\ge a_j^{(m')}$.
 \begin{proof} Because $\rightsquigarrow$ is the transitive closure of $\rightarrow$ it suffices to consider the case when $(i,m)\rightarrow (j,m-1)$ or $(i,m)\rightarrow (j,m+1)$.
 The former case implies that $\mathbf X_{ij}^{(m-1)}\neq 0$ so that, in order for $\lim_{t\rightarrow 0}\lambda(t)\cdot x$ to exist, $a_i^{(m)}-a_j^{(m-1)}$ must be non-negative. The latter case implies that 
 $\mathbf Y_{ij}^{(m)}\neq 0$ so that $a_i^{(m)}-a_j^{(m+1)}$ is non-negative.
  \end{proof}
 \end{lemm}

 Let $x\in \us$ and let $\Lambda_x$ be the set generated by $\{(i,0)\, |\, w_i\neq 0\}$ and the relation $\rightsquigarrow$. That is, 
\[
  \Lambda_x\bydef \{(i,0)\, |\, w_i\neq 0\}\\\cup\left\{\left. (k,m)\in\{1,\ldots,n\}\times \Z_l\, \right|\, (i,0)\rightsquigarrow (k,m)\textrm{ for some $i$ such that }w_i\neq 0\right\}.
 \] Note that $\Lambda_x$ is empty if and only if $\mathbf w=0$.
 
 \begin{lemm}\label{lemm:G.optimal}
  If $\lambda=\left(a_i^{(n)}\right)_{i,n}$ is optimal for $x$ and $(k,m)\in \Lambda_x$ then $a_k^{(m)}=0$.
 \begin{proof}
  Using Corollary \ref{coro:G.zerosandones} and Lemma \ref{lemm:G.rightsquigarrow}, there is an $i$ such that $0\le a_k^{(m)} \le a_i^{(0)} = 0$.
 \end{proof}
 \end{lemm}

 Next, define a one-parameter-subgroup, $\lambda_x$, by
 \[
  \lambda_x\bydef \left( a_k^{(m)}\, \left|\, a_k^{(m)}=\begin{cases} 0&\textrm{ if }(k,m)\in \Lambda_x \\ 1&\textrm{ if }(k,m)\notin \Lambda_x\end{cases}\right.\right).
 \]

 \begin{prop}\label{prop:G.candidatesareoptimal}
  If $x\in \us$ then $\mu(x,\lambda_x)<0$ and $\lambda_x$ is optimal for $x$.
 \begin{proof}
  Let $\left(a_k^{(m)}\right)_{k,m}$ denote the coordinates of $\lambda_x$ as above. Suppose that $\mathbf X_{ij}^{(m)}\neq 0$. Then, whenever $(i,m+1)\in \Lambda_x$, $(j,m)\in \Lambda_x$; therefore,
  $\wt (\lambda_x(t)\cdot x)|_{\mathbf X_{ij}^{(m)}} =  a_i^{(m+1)}-a_j^{(m)}\ge 0$.
  Suppose that $\mathbf Y_{ij}^{(m)}\neq 0$. Then, whenever $(i,m)\in \Lambda_x$, $(j,m+1)\in \Lambda_x$, so
 $  \wt (\lambda_x(t)\cdot x)|_{\mathbf Y_{ij}^{(m)}} =  a_i^{(m)}-a_j^{(m+1)}\ge 0$.
  For all $i$,$\wt (\lambda_x(t)\cdot x)|_{\mathbf v_{i}} =  a_i^{(0)}\ge 0$, and if $w_i\neq 0$ then $(i,0)\in\Lambda_x$ so
 $\wt (\lambda_x(t)\cdot x)|_{\mathbf w_{i}} =  a_i^{(0)}= 0$; therefore, $\lim_{t\rightarrow 0} \lambda_x(t)\cdot x$ exists. 
 
 Suppose that $\mu(x,\lambda_x)=0$. Then $\Lambda_x=\{1,\ldots,n\}\times \Z_l$. Since $x\in \us$, there exists
 some $\nu=\left(b_i^{(m)}\right)_{i,m}\in I$ such that $\mu(x,\nu)<0$. Now, for any $(j,m)$ there exists an $i$ such that
 $(i,0)\rightsquigarrow (j,m)$ and $w_i\neq 0$. By Lemma \ref{lemm:G.rightsquigarrow}, $0=b_i^{(0)}\ge b_j^{(m)}$. This implies that 
 $\nu=\lambda_x$ and so contradicts $\nu \in I$;
 therefore, $\mu(x,\lambda_x)<0$.
 
 Suppose that $\nu=\left(b_k^{(m)}\right)_{k,m}$ is optimal for $x$. Then, by Corollary \ref{coro:G.zerosandones},
 $\nu\in I$ and by Lemma \ref{lemm:G.optimal}, 
 $a_k^{(m)}=0$ implies that $b_k^{(m)}=0$. Let $N_\nu$ be the number of non-zero coordinates of $\nu$ and $N_{\lambda_x}$ the number of 
 non-zero coordinates of $\lambda_x$. Now, since $\nu$ and $\lambda_x$ only have zeros and ones as entries,
 \[
  \frac{\mu(x,\nu)}{\|\nu\|}=-\sqrt{N_\nu}\ge-\sqrt{N_{\lambda_x}}= \frac{\mu(x,\lambda_x)}{\|\lambda_x\|}
 \] so $\lambda_x=\nu$. It follows that $\lambda_x$ is optimal.

 \end{proof}
 \end{prop}
 
 Given an unstable point, $x\in X$, one now has a recipe for producing optimal one-parameter-subgroups, $\lambda_x$; however, not all the $\lambda\in I$ appear in this way. A description of those that do 
 is necessary to decribe the strata of the unstable locus.
 
 Let $\lambda=\left(a_k^{(m)}\right)_{k,m}\in I$ and let 
 \begin{align*}
 i(\lambda)&\bydef \begin{cases}\min\{ m\in \Z_l\, |\, a_k^{(m)}=1\textrm{ for all }k=1,\ldots,n\}&\textrm{ if such a number exists}\\-\infty&\textrm{ otherwise,}\end{cases}\\
 j(\lambda)&\bydef \begin{cases}\max\{ m\in \Z_l\, |\, a_k^{(m)}=1\textrm{ for all }k=1,\ldots,n\}&\textrm{ if such a number exists}\\\infty&\textrm{ otherwise.}\end{cases}
 \end{align*}
 
 \begin{theo}\label{G::theo:essential one-ps} A one-parameter subgroup, $\lambda=\left(a_k^{(m)}\right)_{k,m}\in I$, is optimal for $\theta=(-1,\ldots,-1)$ if and only if it satisfies one of the following (mutually exclusive) conditions.
 \begin{enumerate}[(i)]
  \item For all $k$ and $m$, $a_k^{(m)}=1$.
  \item The number $i(\lambda)$ equals infinity.
  \item Both $i(\lambda)>0$ and $a_k^{(m)}=1$ for all $i(\lambda)\le m \le j(\lambda)$ and all $k=1,\ldots,n$.
 \end{enumerate} 
 \begin{proof} 
 Let $\lambda=\left(a_k^{(m)}\right)_{k,m}$ satisfy one of the three conditions above. A point $x\in X$ will be constructed so that $x\in \us$ and $\lambda=\lambda_x$, the optimal one-parameter
 subgroup for $x$. Define $x$ as follows.
 \begin{align*}
  \mathbf X_{ij}^{(m)}&\bydef\begin{cases} 0&\textrm{ if $a_j^{(m)}=0$ and $a_i^{(m+1)}=1$}\\ 1&\textrm{ otherwise,}\end{cases}&
 \mathbf Y_{ij}^{(m)}&\bydef\begin{cases} 0&\textrm{ if $a_i^{(m)}=1$ and $a_j^{(m+1)}=0$}\\ 1&\textrm{ otherwise,}\end{cases}\\
  v_i&\bydef 1&
  w_i&\bydef\begin{cases}0&\textrm{ if $a_i^{(0)}=1$,}\\1&\textrm{ if $a_i^{(0)}=0$.}\end{cases}
 \end{align*} Now, for all $t\in \CC^* $, $x$ has been defined so that $\lambda(t)\cdot x$ has non-negative weights and $\lambda\in I$ so $\mu(x,\lambda)<0$. Therefore, $x$ lies in the unstable locus.
 
 First, I claim that whenever $\lambda_x$ has a zero entry, the corresponding entry of $\lambda$ must be zero. Second, I claim the converse: whenever an entry of $\lambda$ is zero,
 so must the corresponding entry of $\lambda_x$.
 
 Let $\left(b_k^{(m)}\right)_{k,m}$ denote the entries of $\lambda_x$ and suppose that $b_k^{(m)}=0$ for some $k$ and $m$. Then $(k,m)\in\Lambda_x$ and so there exists an $i$ such that $w_i\neq 0$ and a sequence
 \[
  (i,0)\rightarrow (i_1,m_1)\rightarrow \ldots \rightarrow (i_p,m_p)\rightarrow (k,m).
 \] This implies that $0=a_i^{(0)}$ and, for each $r=1,\ldots,p$, either $X_{i_ri_{r+1}}^{(m_r-1)}\neq 0$ or $Y_{i_ri_{r+1}}^{(m_r)}\neq 0$. Either way, $a_{i_r}^{m_r}\ge a_{i_{r+1}}^{m_{r+1}}$. It follows that
$  0=a_i^{(0)}\ge a_{i_1}^{(m_1)}\ge \cdots \ge a_{i_p}^{(m_p)}\ge a_k^{(m)}$.
 Since $\lambda\in I$, $a_k^{(m)}=0$ and this proves the first claim.
 
 Next suppose that $a_k^{(m)}=0$. Condition (i) doesn't hold now, so $\lambda$ must satisfy one of the other two conditions. This means that $i(\lambda)\neq 0$ and so there must exist an $i$ such that $a_i^{(0)}=0$.
 Suppose that $i(\lambda)$ is infinite. Then for each $m'\in\Z^l$ there exists a $p_{m'}$ such that $a_{p_{m'}}^{(m')}=0$ so that
 $(i,0)\rightarrow (p_1,1) \rightarrow \cdots \rightarrow (p_{m-1},m-1) \rightarrow (k,m)$;
 thus, $(i,0)\rightsquigarrow (k,m)$ and hence $b_k^{(m)}=0$.  
 
 Suppose, instead, that $i(\lambda)=r>0$ and $j(\lambda)=s>0$ are finite. Then Condition (iii) must hold, so every entry between $a_i^{(r)}$ and $a_j^{(s)}$ is $1$ for all $i$ and $j$;
 so either $m<r$ or $m>s$. As in the last case, for 
 each $m'<r$ or $m'>s$, there exists a $p_{m'}$ such that $a_{p_{m'}}^{(m')}=0$. If $m<r$ then, as before,
 $(i,0)\rightarrow (p_1,1) \rightarrow \cdots \rightarrow (p_{m-1},m-1) \rightarrow (k,m)$,
  whereas if $m>s$ then
 $(i,0)\rightarrow (p_{l-1},l-1) \rightarrow \cdots \rightarrow (p_{m+1},m+1) \rightarrow (k,m)$.
  In either case $(i,0)\rightsquigarrow (k,m)$, so $b_k^{(m)}=0$. This proves the second claim. Since $\lambda$ was assumed to be in $I$, the entries of both $\lambda$ and $\lambda_x$ are one when they are non-zero. 
 The two claims therefore show that $\lambda=\lambda_x$. This completes the proof in one direction.
 
 It remains to show that an arbitrary optimal one-parameter subgroup in $I$ must satisfy one of the three conditions. Suppose 
 that some $\lambda\in I$ does not satisfy any of the three conditions, then it suffices to show that $\lambda\neq \lambda_x$ for any $x\in X_\lambda$.
 There are two ways in which both Conditions (ii) and (iii) can fail to hold for $\lambda$. Treat these cases seperately.
 
 First, suppose that $i(\lambda)=0$, so that $a_i^{(0)}=1$ for all $i=1,\ldots,n$. Then any $x\in X_\lambda$ must have $\mathbf w=0$, because the weights of $\lambda$ on non-zero coordinates of $\mathbf w$ would be $-1$. This implies that $\Lambda_x =\varnothing$. However, if $x$ is such that
 $\Lambda_x =\varnothing$ then $\lambda_x=(1,\ldots,1)$ and, since $\lambda$ is not allowed to satisfy Condition (i), $\lambda\neq \lambda_x$.
 
 Second, if $i(\lambda)\neq 0$, then for Condition (iii) to fail there must exist an $m$ such that 
 $1\le i(\lambda)< m < j(\lambda)\le l-1$ and an $i$ such that $a_i^{(m)}=0$. Let $x\in X_\lambda$ be an arbitrary point that is unstable for $\lambda$ and 
 define the entries $(b_i^{(n)})\bydef\lambda_x$. I claim that $\lambda\neq \lambda_x$.
 Suppose that $\lambda=\lambda_x$ so that $b_k^{(m)}=0$. Then there exists some sequence
 $(i,0)\rightarrow (i_1,m_1) \rightarrow \cdots \rightarrow (i_p,m_p)\rightarrow (k,m)$.
 For each $r=1,\ldots,p-1$, $m_{r+1}=m_r \pm 1$, so there must exist a $k'$ such that $m_{k'}=r$ or $m_{k'}=s$.  
 Assume that $m_{k'}=r$, then $(i,0)\rightsquigarrow (i_{k'},r)$, 
 which by Lemma \ref{lemm:G.rightsquigarrow} implies that $0=a_i^{(0)}\ge a_{i_{k'}}^{(r)}=1$. This contradicts the assumption that $b_k^{(m)}=0$, so the claim is proved.
 A similar argument provides a contradiction when there exists a $k'$ such that $m_{k'}=s$. This completes the other direction of the proof.

 \end{proof}
 \end{theo}

 \subsubsection{The Case for other GIT Parameters}
 
The following proposition shows that the optimal one-parameter subgroups for the stability condition $(1,\ldots,1)$ are the same with all the entries multiplied by $-1$.
 \begin{prop}\label{G:prop:negitvies of optimals are optimal} Let $\theta=(\theta_0,\ldots,\theta_{l-1})$ and $\lambda$ a $\theta$-optimal one-parameter subgroup. Then $-\lambda$ is $-\theta$-optimal.
 \begin{proof}
 Suppose that $\lambda$ is $\theta$-optimal. Let $x\in X_\lambda^\theta$, a point that $\lambda$ destabilises with respect to $\theta$ and define a point
 $\hat x =( \hat{\mathbf X}^{(m)}, \hat{\mathbf Y}^{(m)}; \hat v, \hat w)$ by
 \begin{align*}
 \hat{\mathbf X}_{ij}^{(m)}&\bydef \begin{cases} 1 &\textrm{ if }\mathbf Y_{ji}^{(m)}\neq 0 \\ 0&\textrm{ otherwise,}\end{cases}&
 \hat{\mathbf Y}_{ij}^{(m)}&\bydef \begin{cases} 1 &\textrm{ if }\mathbf X_{ji}^{(m)}\neq 0 \\ 0&\textrm{ otherwise,}\end{cases}\\
 \hat v_i&\bydef \begin{cases} 1 &\textrm{ if }\mathbf w_i\neq 0 \\ 0&\textrm{ otherwise,}\end{cases}&
 \hat w_i&\bydef \begin{cases} 1 &\textrm{ if }\mathbf v_i\neq 0 \\ 0&\textrm{ otherwise.}\end{cases}
 \end{align*}
 
 Now, by construction, the weights of any one-parameter subgroup, $\rho$ say, acting on $\hat x$ are precisely the negatives of its weights on $x$ so 
 $\lim_{t\rightarrow 0} \rho(t)\cdot \hat x=
 \lim_{t\rightarrow 0} -\rho(t)\cdot x$, whenever it exists. This means that $\mu^{-\theta}(\hat x,\rho)=\mu^{\theta}(x,-\rho)$; in particular, $\mu^{-\theta}(x,-\lambda)=\mu^\theta(x,\lambda)<0$. Hence $-\lambda$ destabilises $\hat x$, that is $\hat x \in X_{-\lambda}^{-\theta}$.
Suppose that $\nu$ is $-\theta$-optimal for $\hat x$. Then, $x$ is $\theta$-unstable for $-\nu$ so
\begin{align*}
\frac{\mu^{-\theta}(\hat x,-\lambda)}{\|-\lambda\|} = \frac{\mu^{\theta}(x,\lambda)}{\|\lambda\|} \le\frac{\mu^\theta(x,-\nu)}{\|-\nu\|} =\frac{\mu^{-\theta}(\hat x,\nu)}{\|\nu\|}
\end{align*} Since $\nu$ is assumed to be optimal, $\nu=-\lambda$ and the proposition is proved.
 \end{proof}
 \end{prop}
  
\begin{coro}\label{G::prop:negativeoptimals} The set of optimal one-parameter subgroups for $-\theta$ is in bijection with the set of those for $\theta$. In fact,
\[
\Gamma_{-\theta}= \{ -\lambda\, |\, \lambda\in \Gamma_\theta\}
\]
 \begin{proof} Applying the proposition to $\Gamma_\theta$ and $\Gamma_{-\theta}$ gives inclusions of sets in both directions.
 \end{proof}
 \end{coro} 
 
 \begin{defi} Say that a one-parameter subgroup is \textbf{essential} if it is optimal for either $\theta=\pm(1,\ldots,1)$.
 \end{defi}

 \subsection{The McGerty--Nevins Criterion}

Let $\lambda_1,\ldots,\lambda_q$ be the optimal one-parameter subgroups indexing the associated Kirwan--Ness strata of $\us$. For each $k=1,\ldots,q$, let $\mathrm{wt}_{\mathfrak n^-}(\lambda_k)$ denote the
sum of the negative weights of $\lambda_k$ acting on $\mathfrak g$. Let $z\in Z_{\lambda_k}$ and consider the fibre, $N_z$, of the normal bundle of 
$Z_{\lambda_k}$ sitting inside $X$ over $z$. Let $\mathrm{abs.wt}(\lambda_k)$ denote the sum of the absolute values of the $\lambda_k$ weights
on $N_z$. Finally, define $I_k$ to be the set of weights of $\lambda_k$ acting on the symmetric algebra $\mathrm{Sym}^\bullet(N_z)$. These are both independent of the choice of $z\in Z_{\lambda_k}$. 
Let 
\[
\mathrm{shift}(\lambda_k)=- \mathrm{wt}_{\mathfrak n^-}(\lambda_k) - \fr 4 \mathrm{abs.wt}(\lambda_k)\footnote{This is the negative of the shift defined in \cite{MN2}, because their choice of quantised comoment map, $\mu_\mathrm{can}$, is the negative of mine, $\hat \tau$.}.
\]
\begin{theo}(McGerty--Nevins, \cite[Theorem 7.4]{MN2})\label{G::theo:MNT} If, for all $k=1,\ldots,q$,
  \[
\chi(d\lambda_k(1)) \notin \mathrm{shift}(\lambda_k)- I_k
\] then, $\chi$ is good for $\theta$.
\end{theo}

For essential optimals, the condition reduces to checking the following formula. Let $\lambda_\mathbf i$ be an essential optimal one-parameter subgroup so that its differential can be written
 \[
 d\lambda_\mathbf i(1)\bydef \pm(0\cdots0\underbrace{1\cdots1}_{i_0})(0\cdots0\underbrace{1\cdots1}_{i_1})\cdots (0\cdots0\underbrace{1\cdots1}_{i_{l-1}})
 \] for some $\mathbf i=(i_0,\ldots,i_{l-1})\in \Z_{\ge 0}^l$.
 
 \begin{prop}\label{G::prop:Ki's for GL1n} The shifts for essential optimal subgroups are given by the following formula.
 \[
\mathrm{shift}(\lambda_k) =  -\fr 2 i_0 + \sum_{t=0}^{l-1}(n-i_t)(i_t-i_{t-1}).
\] 
 \begin{proof}
Each $\lambda_k$ acts on $X$ with eigenvalues $\{-1,0,1\}$, so $I_k = \Z_{\ge 0}$ and $\mathrm{abs.wt}(\lambda_k)=\dim N_z = \dim X - \dim Z_k$. The action of $\lambda_k$ on $X$ splits each $\mathbf X^{(t)}\in \mathrm{Mat}_{n\times n}(\CC)$ into four blocks of size
$(n-i_t)(n-i_{t+1})$, $(n-i_t)i_{t+1}$, $i_t(n-i_{t+1})$ and $i_ti_{t+1}$ respectively, with each block having the respective weight $0$, $-1$, $1$ and $0$. A count now gives $\mathrm{abs.wt}(\lambda_k)= 2i_0 + 4\sum_t(n-i_t)i_{t+1}$. Similarly,
$\mathrm{wt}_{\mathfrak n^-}(\lambda_k)=-\sum_t(n-i_t)i_t$.
 \end{proof}
 \end{prop}

\section{Bad Parameters for $\mathfrak S_n$, $\mu_3$ and $B_2$}
 Now consider some concrete examples.
 
 \subsection{Localisation for $W=\mathfrak S_n$}
Fix $l=1$. 
\begin{theo}\label{G::theo:bad.params.for. S_n} For $W=\mathfrak S_n$ there are no bad parameters.
\begin{proof} 
  For $\theta<0$ $\{  \lambda_{i}\, |\, i=1,\ldots,n\}$ gives a complete set of representatives for optimal one-parameter subgroups and for $\theta>0$ the optimal one-parameter subgroups are $\{ -\lambda_{i}\}$. For each stratum the weights on the fibres of the normal bundle are $-1$ so the McGerty--Nevins condition reduces to $\chi(d\lambda_i(1))\notin \fr2 i + \Z_{\le 0}$. For $\theta>0$, 
 $\chi(d\lambda_i(1))=-i(c_0+\frac{1}{2})$ so that $c_0$ is bad for $\theta>0$ if $c_0 \in \frac{1}{i}\Z_{\ge 0}$.
  For $\theta <0$, $\chi(d(-\lambda_i)(1))=i(c_0+\frac{1}{2})$ so that $c_0$ is bad if $c_0 \in -1 + \frac{1}{i}\Z_{\le 0}$, see \cite[Corollary 8.1]{MN2}.
  Since these two sets are disjoint there are no bad parameters for this case.
  
  \end{proof}
  \end{theo}
 
 \subsection{Localisation for $W=\mu_3$}\label{Sect:mu3}
 
 Let $n=1$. Define a one-parameter subgroup $\lambda_{i,j}\bydef \lambda_\mathbf k$ where $k_t\bydef \begin{cases} 0 &\textrm{ if $0\le t<i$ or $j\le t<l$}\\ 1&\textrm{ if $i\le t <j$}.\end{cases}$.
 For $\theta=(-1,\ldots,-1)$, the set of one-parameter subgroups is $\{ \lambda_{i,j}\, |\, 1\le i < j \le l\}\cup\{ \lambda_{0,l}\}$. Then
 \[
 \mathrm{shift}(\lambda_{i,j})=\begin{cases} -\fr2 & \textrm{ if $i=0$}\\ -1&\textrm{ if $i\neq 0$}\end{cases}
 \]
 
 Now, in terms of the hyperplane parameters $\{k_1,\ldots,k_{l-1}\}$,
 \[
 \chi(\lambda_{i,j}) = \sum_{t=i}^{j}\left(k_{1-t}-k{-t}+\frac{1}{l}\right)= \frac{j-i+1}{l} + k_{1-i}-k_{-j},
 \] for $i>0$, and
 \[
 \chi(\lambda_{0,l-1}) = \frac{l-1}{l}-k_1.
 \]
 
 \begin{theo}\label{G::theo:bad.params.for.mu3} For $W=\mu_3$ there are no bad parameters.
 \end{theo}
 
 The proof occupies the rest of this section. According to Theorem \ref{G::theo:chi and params}, given an optimal one parameter subgroup $\lambda=(a)(b)(c)$,
 \begin{equation}\label{Equa:mu3}
 \chi(\lambda)= -ck_1 + (c-b)k_2 + \fr 3 (b+c).
 \end{equation} In order for a parameter to be bad for a particular $\theta$, this number must belong to some set of positive real numbers.
 Divide $\R^2$ into the following subsets
 \begin{align*}
 A&\bydef \{ (k_1,k_2)\in \R^2\, |\, k_2\ge \tfrac{4}{3},\, k_2-k_1 \ge \tfrac{2}{3},\, k_1\le \tfrac{7}{6}\} & B&\bydef \{ (k_1,k_2)\in \R^2\, |\, k_2\le -\tfrac{2}{3},\, k_2-k_1 \le -\tfrac{4}{3},\, k_1\ge \fr 6\}\\
C&\bydef \{ (k_1,k_2)\in \R^2\, |\, k_2-k_1 \ge \tfrac{2}{3},\, k_1> \tfrac{7}{6}\} & D&\bydef \{ (k_1,k_2)\in \R^2\, |\, k_2-k_1 \le -\tfrac{4}{3},\, k_1< \fr 6\}.
 \end{align*}
 
 Now, define the following subsets; each is an infinite union of two-dimensional real planes in $\CC^2$.
 \begin{align*}
 B_{1}&\bydef\left\{ (k_1,k_2)\in \CC^2\, \middle|\, k_2-k_1 \in -\tfrac{4}{3} + \Z_{\le 0}\right\} & B_{-1}&\bydef \left\{ (k_1,k_2)\in \CC^2\, \middle|\, k_2-k_1 \in \tfrac{2}{3} +\Z_{\ge 0} \right\}\\
B_{2}&\bydef\left\{ (k_1,k_2)\in \CC^2\, \middle|\, k_2 \in \tfrac{4}{3} + \Z_{\ge 0}\right\}  &B_{-2}&\bydef \left\{ (k_1,k_2)\in \CC^2\, \middle|\, k_2 \in \tfrac{-2}{3} -\Z_{\ge 0} \right\} \\
B_{3}&\bydef\left\{ (k_1,k_2)\in \CC^2\, \middle|\, k_1 \in \tfrac{7}{6} +\fr2 \Z_{\ge 0}\right\} &B_{-3}&\bydef \left\{ (k_1,k_2)\in \CC^2\, \middle|\, k_1 \in \tfrac{1}{6} -\fr2\Z_{\ge 0} \right\}  .
\end{align*} 
\begin{rema}\label{G::rema:compare B_i's with optimals}
These are the sets of bad parameters for the one parameter subgroups as follows.
\begin{align*}
B_1&=B_{(0)(0)(1)} & B_2&=B_{(0)(1)(0)} & B_3&=B_{(0)(1)(1)}\cup B_{(1)(1)(1)}\\
B_{-1}&=B_{(0)(0)(-1)} & B_{-2}&=B_{(0)(-1)(0)} & B_{-3}&=B_{(0)(-1)(-1)}\cup B_{(-1)(-1)(-1)}.
\end{align*} 
\end{rema}

\begin{lemm} There are an inclusions of sets
\begin{align*}
(B_1\cup B_2 \cup B_3)\cap(B_{-1} \cup B_{-2} \cup B_{-3}) &\subset A\cup B\cup C\cup D\\
(B_1\cup B_{-2}\cup B_3)\cap(B_{-1} \cup B_2 \cup B_{-3}) &\subset \R^2 \setminus (A\cup B).
\end{align*}
\begin{proof}
First, note that the sets on the left hand side comprise only points in $\R^2\subset \CC^2$. Indeed, I claim that $B_i\cap B_j\in \R^2$ for any $i\neq j$. 
If $i=\pm 1$ then a point $(k_1,k_2)$ must satisfy, $k_2-k_1 \in \R$, if $i=\pm 2$ then $k_2\in \R$ and if $i=\pm 3$ then $k_1\in \R$. If $j=-i$ then the set is empty. Otherwise, 
for points in $B_i\cap B_j$, two out of three of these conditions must be met and so the claim follows.

Now it is a straight-forward check to see that, 
\[
B_1\cap B_{-2},\, B_1\cap B_{-3},\, B_2\cap B_{-1},\, B_2\cap B_{-3},\, B_3\cap B_{-1},\, B_3\cap B_{-2} \subset A \cup B \cup C\cup D,
\] and a similar check shows the
second inclusion.

\end{proof}
\end{lemm}

\begin{coro}\label{coro:bound}
If a parameter $(k_1,k_2)\in \CC^2$ is bad for every $\theta$ then it belongs to the set $A\cup B \cup C \cup D\subset \R^2$.
\begin{proof}
I claim that if a point does not belong to this set then it is good for either $\theta=\pm(1,1,1)$. Indeed, each of the optimal one-parameter subgroups for these values of $\theta$ 
is accounted for in Remark \ref{G::rema:compare B_i's with optimals} and so $(B_1\cup B_2 \cup B_3)\cap(B_{-1} \cup B_{-2} \cup B_{-3})$ contains all points which are bad for $\theta= \pm(1,1,1)$.
\end{proof}
\end{coro} 
It remains to find, for each parameter $\mathbf k\in A\cup B \cup C \cup D$, some $\theta$ such that $\mathbf k$ is good for $\theta$.
 
 \begin{prop}\label{G::prop:optimals ofr mu3} The optimal one parameter subgroups for $(0,2,-1)$, $(0,-2,1)$, $(0,1,-2)$ and $(0,-1,2)$ are given in Table \ref{G::tabl:optimals for mu3}.
 \begin{proof} See Section \ref{Rema:proofs}.
 \end{proof}
 \end{prop}

\begin{table}[!ht]
 \centering
 \rowcolors{1}{white}{beige1}
 \begin{tabular}{lllllll}
 \toprule
$\theta$& $(-1,-1,-1)$& $(1,1,1)$& $(0,-2,1)$ & $(0,2,-1)$ & $(0,-1,2)$ & $(0,1,-2)$ \\
 \midrule  
&$(0)(0)(1)$ &$(0)(0)(-1)$ &$(0)(0)(-1)$ &$(0)(0)(1)$ &$(0)(0)(-1)$ &$(0)(0)(1)$  \\
			&$(0)(1)(0)$ &$(0)(-1)(0)$ &$(0)(1)(0)$ &$(0)(-1)(0)$ &$(0)(1)(0)$ &$(0)(-1)(0)$  \\
			&$(0)(1)(1)$ &$(0)(-1)(-1)$ &$(0)(1)(1)$ &$(0)(-1)(-1)$ &$(0)(-1)(-1)$ &$(0)(1)(1)$  \\
		&$(1)(1)(1)$ &$(-1)(-1)(-1)$ &$(1)(1)(1)$ &$(-1)(-1)(-1)$ &$(-1)(-1)(-1)$ &$(1)(1)(1)$  \\
			& &				 &$(0)(2)(-1)$ &$(0)(-2)(1)$ &$(0)(1)(-2)$ &$(0)(-1)(2)$  \\
			& &				 &$(-1)(4)(1)$ &$(1)(-4)(-1)$ &$(1)(1)(-4)$ &$(-1)(-1)(4)$  \\
			&  &				 &$(1)(1)(-1)$ &$(-1)(-1)(1)$ &$(-1)(1)(-1)$ &$(1)(-1)(1)$  \\
 \bottomrule \end{tabular}  \caption{ A comparison of the optimal one-parameter subgroups for $W=\mu_3$ and various $\theta$.}\label{G::tabl:optimals for mu3}
 \end{table}
 
Now, I claim that parameters in region $A$ are good for $\theta = (0,1,-2)$, parameters in $B$ are good for $\theta = (0,-1,2)$, parameters in $C$ are good for $\theta = (0,2,-1)$
and parameters in $D$ are good for $\theta = (0,-2,1)$. Indeed, using Equation \ref{Equa:mu3}, one can see that for any of the optimal subgroups, $\lambda$ say, $\chi(d\lambda)>0$ for any parameter $\mathbf k$ in that region.

\subsection{Localisation for $W=G(2,1,2)=B_2$}

This section is a proof of the following theorem.
\begin{theo}\label{G::theo:B_2 params} There are no bad parameters for $G(2,1,2)$.
\end{theo}

Proposition \ref{G::prop:Ki's for GL1n} reduces to the following formula.
 \[
 \mathrm{shift}(\lambda_{i_0,i_1}) = -\fr2 i_0+(n-i_0)(i_0-i_{1})+(n-i_1)(i_1-i_0)=-\fr2 i_0-(i_1-i_0)^2<0.
 \]

Define the following subsets of $\CC^2$.
\begin{align*}
B_1 &\bydef \left\{ \cc\in \CC^2\, |\, c_1\in -\fr2 + \fr 2\Z_{\le 0}\right\}\setminus\{c_1=-1\} & B_{-1} &\bydef \left\{ \cc\in \CC^2\, |\, c_1\in \tfrac{3}{2}+ \fr2 \Z_{\ge 0}\right\}\setminus\{c_1=2\} \\
B_2 &\bydef \left\{ \cc\in \CC^2\, |\, c_0\in   -1+\fr2\Z_{\le 0}\right\} & B_{-2} &\bydef \left\{ \cc\in \CC^2\, |\, c_0\in  \fr2\Z_{\ge 0}\right\} \\
B_3 &\bydef \left\{ \cc\in \CC^2\, |\, c_0-c_1\in -\tfrac{5}{2}+ \Z_{\le 0}\right\} & B_{-3} &\bydef \left\{ \cc\in \CC^2\, |\, c_0-c_1\in \tfrac{1}{2}+ \Z_{\ge 0}\right\} \\
B_4 &\bydef \left\{ \cc\in \CC^2\, |\, c_0+c_1\in -\tfrac{3}{2} + \Z_{\le 0}\right\} & B_{-4} &\bydef \left\{ \cc\in \CC^2\, |\, c_0+c_1\in \tfrac{3}{2}+ \Z_{\ge 0}\right\} 
\end{align*}

Given a one-parameter subgroup $\lambda$, let $B_{d\lambda(1)}$ be the set of parameters that are bad for $\lambda$. Here,
\begin{align*}
B_1 &=B_{(00)(01)}\cup B_{(00)(11)}, &
B_{-1} &= B_{(00)(0-1)}\cup B_{(00)(-1-1)},\\ 
B_2 &=B_{(01)(01)}\cup B_{(11)(11)}, &
B_{-2} &=B_{(0-1)(0-1)}\cup B_{(-1-1)(-1-1)},\\ 
 B_3 &=B_{(01)(00)}, &
 B_{-3} &=B_{(0-1)(00)},\\ 
B_4 &=B_{(01)(11)},&
B_{-4} &=B_{(0-1)(-1-1)},
\end{align*} so that $B_1\cup\cdots\cup B_4$ and $B_{-1}\cup\cdots\cup B_{-4}$ are the sets of bad points for the GIT parameters $(-1,-1)$ and $(1,1)$ respectively. 

Define the following disjoint subsets of $\R^2$, thought of as sitting inside the copy of $\CC^2$ above.
\begin{align*}
A&\bydef \left\{ \cc\in \R^2\, \middle|\, c_0>-1,\, c_1-c_0\ge \tfrac{5}{2}\right\} &
B_+&\bydef \left\{ \cc\in \R^2\, \middle|\, c_0\le -1,\, c_0+c_1\ge \tfrac{1}{2}\right\}\\
B_-&\bydef \left\{ \cc\in \R^2\, \middle|\, c_1 \ge \tfrac{3}{2}, \, c_1-c_0>\tfrac{1}{2}\right\} &
C&\bydef \left\{ \cc\in \R^2\, \middle|\, c_0<0,\, c_0-c_1\ge \tfrac{1}{2}\right\} \\
D_+&\bydef \left\{ \cc\in \R^2\, \middle|\, c_1\le -\fr2,\, c_0+c_1\ge -\fr2 \right\} &
D_-&\bydef \left\{ \cc\in \R^2\, \middle|\, c_0\ge 0,\, \,c_0-c_1 > \tfrac{1}{2}\right\}.
\end{align*} Let $B=B_+\cup B_-$ and $D=D_+\cup D_-$.

\begin{lemm}
The set, $A\cup B\cup C\cup D$, gives a bound for the set of bad parameters when $W=G(2,1,2)$.
\begin{proof}
The proof follows the same argument in Corollary \ref{coro:bound}.
\end{proof}
\end{lemm}

By the lemma above, it suffices to find, for each region, $A,\ldots,D_-$, some $\theta$ such that parameters in that region are good for that $\theta$.

\begin{prop}\label{G::theo:optimals for B_2} The optimal one-parameter subgroups for $\theta=\pm(1,-2)$, $\pm(3,-2)$, $\pm(3,-1)$ are given in Table
\ref{G::tabl:optimals for B_2}.
\begin{proof} See Section \ref{Rema:proofs}
\end{proof}
\end{prop}

\begin{table}[!ht]
 \centering
 \rowcolors{1}{white}{beige1}
 \begin{tabular}{lllllllll}
 \toprule
$\theta$&$(-1,-1)$&$(1,1)$& $(-1,2)$& $(1,-2)$& $(-3,2)$ & $(3,-2)$ & $(-3,1)$ & $(3,-1)$ \\
 \midrule  
&(00)(01)  		&(00)({-1}0)			&(00)({-1}0)		&(00)(01)  		&(00)({-1}0)			&(00)(01)			&(00)({-1}0)			&(00)(01)\\
&(01)(00)		&(-10)(00)				&(01)(00)			&(-10)(00)		&(01)(00)				&(-10)(00)			&(01)(00)				&(-10)(00)\\
&(01)(01)		&({-1}0)({-1}0)			&({-1}0)({-1}0)		&(01)(01)		&(01)(01)				&({-1}0)({-1}0)		&(01)(01)				&({-1}0)({-1}0)\\
&(00)(11)   	&(00)({-1}{-1})			&(00)({-1}{-1})		&(00)(11)   	&(00)({-1}{-1})			&(00)(11)			&(00)({-1}{-1})			&(00)(11)  \\
&(01)(11)		&({-1}0)(-1{-1})		&({-1}0)(-1{-1})	&(01)(11)		&({-1}0)(-1{-1})		&(01)(11)			&(01)(11)				&({-1}0)(-1{-1})\\
&(11)(11)		&(-1{-1})(-1{-1})		&(-1{-1})(-1{-1})	&(11)(11)  		&(11)(11)				&(-1{-1})(-1{-1})	&(11)(11)				&(-1{-1})(-1{-1})\\
&&													&(11)(-2{-2}) 	&(-1{-1})(22) 		&(33)(-2{-2})			&(-3{-3})(22)		&(33)(-1{-1})			&(-3{-3})(11)\\
&&										&(01)(-2{-2})   &({-1}0)(22)		&(03)(-2{-2})			&(-30)(22)			&(03)(-1{-1})			&(-30)(11)\\
&&										&(-12)(-4{-1}) 	&(-21)(14)			&(16)(-41)				&(-6{-1})(-14)		&(03)(-10)				&(-30)(01)\\
&&										&(01)(-20)		&(-10)(02)   		&(03)(-20)				&(-30)(02)			&(13)(-11)				&(-3{-1})(-11)\\
&&										&(-10)(-4{-1})	&(01)(14) 			&(22)(-32)				&(-2{-2})(-23)		&(19)(11)				&(-9{-1})(-1{-1})\\
&&										&(-11)(-1{-1}) 	&(-11)(11) 			&(-19)(-1{-1})			&(-91)(11)			&(55)(-35)				&(-5{-5})(-53)\\
&&			&&						&(01)(-41)				&(-10)(-14)			&(01)(-11)				&(-10)(-11)\\
 \bottomrule \end{tabular}\label{G::tabl:optimals for B_2}  \caption{ A comparison of the optimal one-parameter subgroups for $W=B_2$ and various $\theta$.}
 \end{table}

By Theorem \ref{G::theo:chi and params}, 
\begin{equation}\label{Equa:B2}
\chi((ab)(cd))= (a+b)c_0 + (c+d-a-b)c_1 +(a+b- \fr2 (c+d)).
\end{equation}

I claim that parameters in region $A$ are good for $\theta = (1,-2)$, parameters in $B_+$ are good for $\theta = (3,-2)$, parameters in $B_-$ are good for $\theta = (3,-1)$, 
parameters in $C$ are good for $\theta = (-1,2)$, parameters in $D_+$ are good for $\theta = (-3,1)$
and parameters in $D_-$ are good for $\theta = (-3,2)$. Indeed, using Equation \ref{Equa:B2}, one can see that for any of the optimal 
subgroups, $\lambda$ say, $\chi(d\lambda)>0$ for any parameter $\cc$ in that region.

\subsection{The proofs of Propositions \ref{G::prop:optimals ofr mu3} and \ref{G::theo:optimals for B_2}}\label{Rema:proofs}

For the proofs of the two propositions, the reader is referred to Appendix C of \cite{JenThesis}. Candidate optimal subgroups were 
first calculated using a computer program. They are proved to be optimal by the following argument. The set of unstable points are 
classified by `types' based on the criterion of King in the following way. Each unstable point, $p\in \us$, corresponds to a representation 
of $Q^l_\infty$ that has a destabilising subrepresentation with a certain dimension vector; this is referred to as a \textbf{type} 
of the point $p\in \us$. Note that a point may have several different types.
For a point $q\in X$, having one type (or not having another) can put various restrictions on the coordinates of that point. 

The possible types of each unstable point are divided into different cases, each corresponding to a proposed optimal candidate. It is 
shown that the various cases are disjoint and exhaustive of all unstable points. The following lemma is then used to show that all points 
belonging to a particular case have the proposed subgroup as an optimal.

Let $T$ be the subgroup of diagonal matrices in $G$. 
\begin{lemm} Let $q\in \us$ be a point that is specified by its types. Suppose that $\lambda=(a_k^{(m)})\in \mathbb Y(T)$ is an optimal 
one-parameter subgroup for some point $p$ in the orbit of $q$. Then a one-parameter subgroup $\mu$ is optimal for $q$ if the following conditions are satisfied.
\begin{enumerate}
\item From the description of $q$ in terms of types, a restriction is imposed on $\lambda$ that bounds the value of $m(\lambda)$ above.
\item The candidate subgroup $\mu$ has $m(\mu)$ greater than or equal to this bound.
\item The subgroup $\mu$ destabilises some point $r$ in the $G$-orbit of $q$.
\end{enumerate}
\begin{proof} First, note that if $\lambda'\in \mathbb Y(G)$ is optimal for some point $q$ then, because the image of any one-parameter 
subgroup lies in some maximal torus of $G$ and $T$ is a maximal torus, there is some $\lambda\in \mathbb Y(T)$ that is optimal for some 
point $p$ in the same orbit as $q$. The function $M$ is constant on $G$-orbits so the first two conditions imply that $M(p)=m(\lambda) \le m(\mu)$. 
If $\mu$ destabilises a point $r$ in the $G$-orbit of $q$ then $m(\mu)\le M(r)$. Thus $M(r)=m(\mu)$. Thus $r\in S_\mu$. Now, since 
$S_{\langle \mu \rangle}^+=G\times_{P(\mu)} S_\mu$ it follows that $q\in S_{\langle \mu \rangle}^+$.
\end{proof}
\end{lemm}

Here's an example. If $W=\mu_3$ and $\theta=(0,-2,1)$ then every unstable point must correspond to a representation of $\overline{Q^3_\infty}$ 
 with a destabilising subrepresentation with one of the following dimension vectors: $(1,1,1,0)$,
 $(1,0,1,0)$, $(0,1,1,1)$, $(0,1,1,0)$, $(0,0,1,1)$, $(0,0,1,0)$. For the third case, points are described as those that have type
$(0,0,1,0)$, but neither $(1,0,1,0)$ nor $(1,1,1,0)$. Suppose that $q=(X^{(0)},X^{(1)},X^{(2)},Y^{(0)},Y^{(1)},Y^{(2)};v,w)\in X\cong\CC^8$ belongs to this case. Let 
$\lambda=(a)(b)(c)\in \mathbb Y(T)$ be an optimal for $q$. The weights of $\lambda$ on $q$ are $(b-a,c-b,a-c,a-b,b-c,c-d;a,-a)$. Being of type $(0,0,1,0)$ but not being of type $(1,1,1,0)$ implies that
$Y^{(2)}\neq 0$ so that $a\le c$. But also not being of type $(1,0,1,0)$ implies that $v\neq 0$ so that $0\le a \le c$.
Now,
\begin{align*}
0&\le 4a^2+5c^2+2c(2b-c) \\
&= 4(a^2+b^2+c^2)-(2b-c)^2,\\
&= 4\|\lambda\|^2 - \mu(p,\lambda)^2.
\end{align*}
Thus $m(\lambda)\le \sqrt{4} = m((0)(1)(0))$. 
Being of type $(0,0,1,0)$ implies $X^{(1)}=Y^{(0)}=0$ so $(0)(1)(0)$ acts by non-zero weights on $q$ and so is optimal for $q$.

\bibliographystyle{alpha}
\bibliography{bib}
\addcontentsline{toc}{section}{Bibliography}

\end{document}